%% file: thesis.tex
\documentclass[12pt]{amsbook} 
\usepackage{amssymb}
\usepackage{latexsym}
\usepackage{amsfonts}
\usepackage[all]{xy}
\input{commands}

\begin{document}
%front matter
\title{Symmetric Graphs and their Quotients}
\author{Robin Langer}

\maketitle
\tableofcontents
%main matter
\include{introduction}

\include{background}

\include{coset}
\include{geometric}

\include{threearc}

%\include{symmetric-design}

%back matter
%bibliography
\nocite{*}
\bibliographystyle{plain} % or amsalpha (or others)
\bibliography{references}
\end{document}

%% file: commands.tex
\theoremstyle{plain} 
\newtheorem{theorem}{Theorem}
\newtheorem{lemma}{Lemma}

\newtheorem{proposition}{Proposition}

\theoremstyle{definition} 
\newtheorem{definition}{Definition}
\newtheorem{example}{Example}
\newtheorem{question}{Question}

\newtheorem*{condition}{Condition (PE)}
\theoremstyle{remark}

\DeclareMathOperator{\Aut}{Aut}
\DeclareMathOperator{\Cos}{Cos}
\DeclareMathOperator{\Core}{Core}
\DeclareMathOperator{\orb}{Orb}
\DeclareMathOperator{\Cay}{Cay}
\DeclareMathOperator{\Sab}{Sab}
\DeclareMathOperator{\sab}{Sab}

\DeclareMathOperator{\Arc}{Arc}
\DeclareMathOperator{\val}{val}

\DeclareMathOperator{\sub}{Sub}

\newcommand{\cB}{\mathcal{B}}
\newcommand{\cS}{\mathcal{S}}
\newcommand{\cP}{\mathcal{P}}
\newcommand{\cZ}{\mathcal{Z}}
\newcommand{\cD}{\mathcal{D}}

%% file: introduction.tex
\chapter*{Introduction}

% Introductory paragraph:

A graph is a combinatorial object that captures abstractly 
the idea of a relationship amongst the elements in a set.
Associated with every combinatorial object is a group of symmetries, 
or an \emph{automorphism group}.
An automorphism is, loosly speaking, a structure preserving map from the object to itself.
A homomorphism maps a complex object onto a simpler one
in such a way that certain features of the original object are preserved
while others are lost.

In this paper we study a family of highly symmetric graphs
using a mixture of group theoretic and combinatorial techniques.
The graphs we study have the property that locally
they ``look the same'' at every vertex, 
while globally they are rich in structure.
In particular we look at homomorphic images, 
or quotients, of symmetric graphs. 
We would like to understand how the combinatorial structure of a
symmetric graph is related to that of its quotients.

When passing from a graph to its quotient, information is lost.
For any given symmetric graph there are, in fact, infinitely many
larger symmetric graphs which admit the given graph as a quotient.
Where is this information being lost to?
How is it possible to take a symmetric graph and ``unfold''
it into a larger symmetric graph which
admits the original as a quotient?
What extra information is needed?
\hfill \break

In Chapter 1, we introduce the basic notions from the theory of
permutation groups necessary to give the definition of
a symmetric graph, an imprimitive symmetric graph
and the quotient of an imprimitive symmetric graph.

In Chapter 2 we introduce the idea of coset spaces and coset graphs. 
We see that, in some sense, 
symmetric graphs capture combinatorially
the way that a subgroup sits inside a larger group.

In Chapter 3 we look at the ``extension problem''
for symmetric graphs and describe a ``geometric approach''
to the problem suggested by Gardiner and Praeger.

In Chapter 4 we look at a number of methods
for constructing symmetric graphs with a given quotient

%Butterflies in bloom. See how the grow from a small seed into such majestic beings?

%% file: background.tex
\chapter{$G$-Sets and $G$-Graphs}

\section{Permutation Groups}

% Definition of permutation group:
\begin{definition}[Permutation group]
A \emph{permutation group} is a triple $(G, \Omega, \rho)$
where $G$ is a group, $\Omega$ is a set and $\rho$
is a homomorphism:
$$ \rho : G \to \Aut(\Omega). $$
We say that $G$ {\it acts on} $\Omega$ as a group of permutations.
The action is said to be \emph{faithful} if $\ker(\rho) = \{1\}$.
\end{definition}

% Remark on notation:

We will usually neglect to mention the homomorphism $\rho$ explicitely,
and speak of the permutation group $(G, \Omega)$.
We write $\alpha^{\rho(g)}$ or just $\alpha^g$ to indicate
the action of a permutation $g \in G$ on a point $\alpha \in \Omega$. 
This allows us to compose permutations as 
$\alpha^{(gh)} = (\alpha^{g})^{h}$
rather than $(gh)(\alpha) = h(g(\alpha))$.

If $\Delta$ is a subset of $\Omega$ and $g$ is an element of $G$, then we write
$\Delta^g = \{ \alpha^g : \alpha \in \Delta \}$
to indicate the image of $\Delta$ under the action of $g$.

\hfill \break
Sometimes we will call $(G,\Omega)$ a \emph{representation} 
of $G$ as a group of permutations of $\Omega$, rather than a permutation group.
If $|\Omega| = n$ then we say that $(G,\Omega)$ is a permutation
representation of $G$ of degree n.

% Definition: permutation equivalence:
\begin{definition}[Permutation Equivalence]
Two permutation groups $(G_1, \Omega_1)$ and
$(G_2, \Omega_2)$ are said to be \emph{permutation equivalent}
if there is an isomorphism $\varphi : G_1 \to G_2$
and a bijection $\eta : \Omega_1 \to \Omega_2$ 
such that for every $g \in G_1$ the following diagram commutes:

% Commutative diagram using xypic:
\centerline{
 \xymatrix{
 \ar[d]_{\eta} \ar[r]^g \Omega_1 & \Omega_1 \ar[d]_{\eta} \\
 \ar[r]_{\varphi(g)} \Omega_2 & \Omega_2
 }
}
\end{definition} 

Permutation equivalence is an equivalence relation
on the set of permutation representations of given group.
If we relax the condition that $\eta$ is bijective 
we obtain a {\it permutation homomorphism}.
%For a fixed group $G$,
%the set of permutation representations of $G$, 
%together with the permutation homomorphisms form a category,
%which we denote by $G$-Set.

% Definition of transitivity
\begin{definition}[Transitive Permutation Group]
A permutation group  $(G, \Omega)$ is said to be \emph{transitive}
if for any $\alpha, \beta \in \Omega$ there exist a $ g \in G $
such that $\alpha^g = \beta$
\end{definition}

% Definition of orbit
\begin{definition}[Orbit]
If $(G, \Omega)$ is a permutation group and $\alpha \in \Omega$ 
then the \emph{orbit of $\alpha$ under the action of $G$} is:
$$ \alpha^G = \{ \alpha^x : x \in G \} $$
\end{definition}

For a transitive permutation group $\alpha^G = \Omega$.
Otherwise $\Omega$ maybe be partitioned into orbits:
$ \Omega = \coprod_{i = 1..n} \Omega_i $
such that for each $i$, the permutation group $(G,{\Omega_i})$ is transitive.

% category theory
%\hfill \break
%In category theoretic language, the transitive representations
%of a group $G$ are precisely the objects of $G$-Set
%which have no proper subobjects.

\section{Graphs}

The most general definition of a {\it graph} is a pair of sets $(V, A)$
with $A \subseteq V \times V$.
The elements of $V$ are refered to as the {\it vertices} of the graph
and the elements of $A$ are refered to as {\it arcs}.
A {\it simple graph} is a graph such that:
$$ (\alpha,\beta) \in A \Longleftrightarrow (\beta, \alpha) \in A
\mbox{ for all } \alpha, \beta \in V.$$ 

A graph that is not simple is said to be {\it directed}.
The term ``digraph'' is also frequenty used.
An {\it edge} of a graph is an unordered pair $\{\alpha,\beta\}$
where $(\alpha,\beta)$ is an arc and $\alpha \neq \beta$.
A {\it loop} is an arc of the form $(\alpha,\alpha)$.
In this paper we will be mostly interested
in simple graphs without loops
and shall use the term ``graph'' to mean
``simple graph without loops''. We give a formal definition for ease of reference.

\begin{definition} [Graph]
A {\it graph} is a pair of sets $(V,A)$ such that:
$$A \subseteq (V \times V) \diagdown \{ (\alpha,\alpha) : \alpha \in V \}$$

and: 
$$(\alpha,\beta) \in A \Longleftrightarrow (\beta,\alpha) \in A
\mbox{ for all } \alpha, \beta \in V.$$
\end{definition}

When more than one graph is being discussed, we shall write
$\Gamma = (V \Gamma, A \Gamma)$ or $\Sigma = (V \Sigma, A \Sigma)$
so that it is clear to which graph each set of vertices and arcs belongs.
If $\alpha$ is a vertex of some graph $\Gamma$,
then the {\it neighbourhood} of $\alpha$ in $\Gamma$
is the set $\Gamma(\alpha) = \{ \beta \in V \Gamma : (\alpha,\beta) \in A \Gamma \}$.
If $\alpha$ is a vertex of the graph $\Sigma$
then the neighbourhood of $\alpha$ in $\Sigma$ is denoted by $\Sigma(\alpha)$.
We write $E \Gamma$ (or $E \Sigma$) to denote the set of edge of $\Gamma$ (or $\Sigma$)
and $N \Gamma$ (or $N \Sigma$) to denote the set of neighbourhoods.

\begin{definition} [Automorphism of a Graph]
An automorphism of a graph $\Gamma$ is a bijection
$\varphi : V \Gamma \to V \Gamma$ such that:
$$ (\alpha,\beta) \in A \Gamma \Longleftrightarrow
(\varphi(\alpha),\varphi(\beta)) \in A \Gamma \mbox{ for all } \alpha,\beta \in V \Gamma.$$ 
\end{definition}

The autmorphisms of a graph form a group which we denote by $\Aut(\Gamma)$.
If $\rho : G \to \Aut(\Gamma)$ is a homomorphism
then $(G,V \Gamma,\rho)$ is a permutation group
and we say that $G$ acts on $\Gamma$ as a group of automorphisms.

\begin{definition} [Vertex-transitive Graph]
A vertex transitive graph is a triple $(G,\Gamma,\rho)$
where $\rho : G \to \Aut(\Gamma)$ is a homomorphism
and $(G,\Gamma,\rho)$ is transitive.
\end{definition}

As with permutation groups, we will often neglect to mention $\rho$
explicitely. We shall also sometimes say that $\Gamma$ is a $G$-vertex
transitive graph rather than $(G,\Gamma)$ is a vertex transitive graph.
If we say simply that $\Gamma$ is a vertex transitive graph
then we mean that $(\Aut(\Gamma),\Gamma)$ 
is an vertex transitive graph.
%Need to define the valency of a graph.
%Also complete graph and bipartite graph.

\section{Local Properties}

\begin{definition}[Stabilizer]
Let $(G, \Omega)$ be a transitive permutation group,
and let $\alpha$ be any point of $\Omega$.
The \emph{point stabilizer} of $\alpha$ is:
$$ G_a = \{ x \in G : \alpha^x = \alpha  \} $$
\end{definition}

\begin{lemma}
$G_{\alpha^g} = g^{-1} G_{\alpha} g$
for any $g \in G$. 
\begin{proof}
Suppose that $x \in G$ is such that $(\alpha^g)^x = \alpha^g$.
It follows that $\alpha^{gxg^{-1}} = \alpha$
and so $gxg^{-1} \in G_\alpha$. That is $x \in g^{-1}G_\alpha g$.
Conversely suppose that $x \in g^{-1}G_\alpha g$.
Then $gxg^{-1} \in G_\alpha$ and so $\alpha^{gxg^{-1}} = \alpha$.
That is $\alpha^{gx} = \alpha$, so $x \in G_{\alpha^g}$.
\end{proof}
\end{lemma}

Since $(G,\Omega)$ is transitive, every element of $\Omega$
may be expressed in the form $\alpha^g$ for some $g \in G$.
Thus lemma 1 tells us that the point stabilizers of a transitive
representation of $G$ form a family of conjugate subgroups of $G$.

If $G$ acts on $\Gamma$ as a group of automorphisms
and $G_\alpha$ is the stabilizer of the vertex $\alpha$
then $(G_\alpha, \Gamma(\alpha))$ is a permutation group.

\begin{proposition}
if $(G, \Gamma)$ is a vertex transitive graph
then for all $\alpha, \beta \in V \Gamma$
the permutation groups $(G_\alpha, \Gamma(\alpha))$ and
$(G_\beta, \Gamma(\beta))$ are permutation equivalent.

\begin{proof}
Suppose that $\beta = \alpha^g$. 
Let $\eta : \Gamma(\alpha) \to \Gamma(\beta) $ be given by: $ \gamma \mapsto \gamma^g $,
and let $\varphi : G_\alpha \to G_\beta $
be given by: $ x \mapsto g^{-1} x g $
then:
\begin{eqnarray*}
\eta(\gamma^x) & = & \gamma^{xg} \\
               & = & \gamma^{gg^{-1}xg} \\
	       & = & (\gamma^g)^{\varphi(x)} \\
               & = & \eta(\gamma)^{\varphi(x)}
\end{eqnarray*}	       
for any $\gamma \in \Gamma(\alpha)$ and any $x \in G_\alpha$.
Thus the pair of maps $(\eta,\varphi)$ establish the desired
permutation equivalence.
\end{proof}
\end{proposition}

If $\Gamma$ is $G$-vertex transitive,
then for any permutation group theoretic property $\cP$
we say that $\Gamma$ is {\it $G$-locally $\cP$}
if $(G_\alpha, \Gamma(\alpha))$ has property $\cP$.

\section{Symmetric Graphs}

\begin{definition} [Symmetric Graph]

A symmetric graph is a triple $(G,\Gamma,\rho)$
where $G$ is a group, $\Gamma$ is a graph
and $\rho$ is a homomorphism:
$$\rho : G \to \Aut(\Gamma)$$
such that $\Gamma$ is $G$-vertex transitive and $G$-locally transitive.
\end{definition}

Again, we won't always mention $\rho$ explicitely
and will speak of the symmetric graph $(G,\Gamma)$, or sometimes
the $G$-symmetric graph $\Gamma$.

\begin{definition} [s-arc]
An $s$-arc of a graph $\Gamma$ is a sequence of vertices $(v_0, v_1, ... v_s)$
such that $v_i$ is adjacent to $v_{i+1}$ for each $i$ and $v_{i-1} \neq v_{i+1}$.
The set of $s$-arcs of $\Gamma$ is denoted by $\Arc_s(\Gamma)$.
\end{definition}

If $G$ acts on $\Gamma$ as a group of automorphisms,
then $G$ also acts on $\Arc_s(\Gamma)$ in a natural way:
$$ (v_0, v_1, ... v_s)^g = (v_0^g, v_1^g, ... v_s^g)$$
The graph $\Gamma$ is said to be $(G,s)$-arc transitive if $(G, \Arc_s(\Gamma)$
is transitive.
Historically there has been much interest in highly arc-transitive graphs,
that is graphs that are $s$-arc transitive for large $s$.

One of the oldest results in the area is the theorem by Tutte \cite{TUTTE1,TUTTE2},
proved using combinatorial arguments,
that there are no $s$-arc transitive graphs
of valency $3$ for $s > 5$.
More recently Weiss \cite{WEISS} was able to show, with the aid
of the classification of finite simple groups
that there are no $s$-arc transitive graphs,
of any valency, for $s > 7$.

Symmetric graphs have historically been characterized by their arc transitivity
rather than their local transitivity.
A $(G,0)$-arc transitive graph is just a $G$-vertex transitive graph.
The next theorem shows that $(G,1)$-arc transitive graphs
wihtout isolated vertices are symmetric graphs.

\begin{proposition}
if $(G, \Gamma)$ is a symmetric graph, the $\Gamma$ is $(G,1)$-arc transitive.
Conversly, if $\Gamma$ contains no isolated vertices and is $(G,1)$-arc transitive,
then $(G, \Gamma)$ is a symmetric graph.

\begin{proof}
Suppose that $(G, \Gamma)$ is a symmetric graph, then
$G$ acts on $A \Gamma$ in the obvious way:
$(\alpha,\beta)^x = (\alpha^x, \beta^x)$.
Let $(\alpha_1, \beta_1)$ and $(\alpha_2, \beta_2)$
be any two arcs.
By the $G$-vertex transitivity of $\Gamma$ we can find a $g \in G$ such that
$\alpha_1^g = \alpha_2$.
Since $\beta_1 \in \Gamma(\alpha_1)$ and $g$ is an automorphism,
we must have $\beta_1^g \in \Gamma(\alpha_1^g) = \Gamma(\alpha_2)$.
By $G$-local transitivity we can find $h \in G_{\alpha_2}$
such that $(\beta_1^g)^h = \beta_2$.
It follows that
$(\alpha_1, \beta_2)^{gh} = (\alpha_1^{gh}, \beta_1^{gh})
= (\alpha_2, \beta_2)$.
Thus $G$ acts transitively on the arcs of $\Gamma$.

Now suppose that $\Gamma$ has no isolated vertices and $(G,A \Gamma)$
is transitive. For any two vertices $\alpha_1$ and $\alpha_2$,
there is some arc beginning at $\alpha_1$, say $(\alpha_1,\beta_1)$
and some arc beginning at $\alpha_2$, say $(\alpha_2, \beta_2)$.
By arc transitivity we can find some $g \in G$ carrying
$(\alpha_1, \beta_1)$ to $(\alpha_2, \beta_2)$,
and this $g$ must carry $\alpha_1$ to $\alpha_2$.
So $\Gamma$ is $G$-vertex transitive.

Suppose $\gamma_1$ and $\gamma_2$ are both elements of $\Gamma(\alpha)$.
The $(\alpha,\gamma_1), (\alpha,\gamma_2) \in A \Gamma$.
By arc transitivity we can find a $g \in G$ such that
$(\alpha, \gamma_1)^g = (\alpha, \gamma_2)$
and this $g$ must carry $\gamma_1$ to $\gamma_2$.
Thus $\Gamma$ is $G$-locally transitive.
The result follows.
\end{proof}
\end{proposition}

\section{Orbital Graphs}

If $(G,\Omega)$ is a transitive permutation group,
then $G$ has a natural action on $\Omega \times \Omega$
given by:
$$ (\alpha,\beta)^g = (\alpha^g, \beta^g) $$
Although $(G,\Omega)$ is transitive, $(G,\Omega \times \Omega)$
need not be. 
The orbits of $G$ on $\Omega \times \Omega$ are called {\it orbitals}

\begin{lemma} 
If $(G,\Omega)$ is a transitive permutation group
and $\alpha \in \Omega$
then there is a natural correspondence between the orbitals of $(G,\Omega)$
and the orbits of $(G_\alpha,\Omega)$ 
\begin{proof}
Let $\Delta$ be any orbital of $(G,\Omega)$
and let $(\gamma,\delta)$ be any element of $\Delta$.
Let $g \in G$ be such that $\gamma^g = \alpha$.
Since $\Delta$ is an orbital, it follows that
$(\alpha,\delta^g) \in \Delta$.
Thus every orbital of $(G,\Omega)$ contains an element
of the form $(\alpha,\beta)$.

Let $\eta : \Omega \to \Omega \times \Omega$ be the map given by
$\beta \mapsto (\alpha, \beta)$. 
Suppose that $\beta_1$ and $\beta_2$
lie in the same orbit of $(G_\alpha,\Omega)$.
Then there is some $g \in G_\alpha$ such that $\beta_1^g = \beta_2$.
It follows that $(\alpha,\beta_1)^g = (\alpha,\beta_2)$
and so $(\alpha,\beta_1)$ and $(\alpha,\beta_2)$ lie in the same
orbital of $(G,\Omega)$.

Conversely suppose that $(\alpha,\beta_1)$ and $(\alpha,\beta_2)$
lie in the same orbital of $(G,\Omega)$.
Then there is some $g \in G$ such that $(\alpha,\beta_1)^g = (\alpha,\beta_2)$.
That is $\alpha^g = \alpha$ and $\beta_1^g = \beta_2$.
It follows that $g \in G_\alpha$ and $\beta_1$ and $\beta_2$
lie in the same orbit of $(G_\alpha, \Omega)$.
Thus, the map $\eta$ in fact induces a map from the orbits of $(G_\alpha,\Omega)$
to the orbitals of $(G,\Omega)$.
\end{proof}
\end{lemma}

The orbital $\{ (\alpha,\alpha) : \alpha \in \Omega)$
is given a special name. Its called the {\it diagonal} orbital.
The {\it rank} of a permutation group is the number of orbitals.
An orbital $\Delta$ is said to be {\it self-paired} if
$$ (\alpha, \beta) \in \Delta \Longleftrightarrow
   (\beta, \alpha) \in \Delta $$

\begin{definition} [Orbital Graph]
If $(G,\Omega)$ is a transitive permutation group and
$\Delta$ is a self-paired orbital 
then the {\it orbital graph} $\orb_\Delta(G, \Omega)$
is the graph with vertex set $\Omega$ and arc-set $\Delta$
Need to say more about the diagonal orbital.
\end{definition}

Note that if $\Delta$ is taken to be the diagonal orbit,
then the resulting graph has loops, and so is not actually
a graph by our definition. Generally we assume that $\Delta$
is not the diagonal orbit.

\begin{proposition}
Orbital graphs are symmetric. Every symmetric graph is an orbital graph
\begin{proof}
Suppose that $\Gamma = \orb_\Delta(G,\Omega)$ 
for some transitive permutation group
$(G,\Omega)$ and some self-paired orbital $\Delta$.
Since $V \Gamma = \Omega$ and $(G,\Omega)$ is transitive,
$\Gamma$ is $G$-vertex transitive.
Since $A \Gamma = \Delta$ and $\Delta$ is an orbital
$\Gamma$ is $G$-arc transitive.
It follows that $\Gamma$ is symmetric.

Now suppose that $\Gamma$ is a $G$-symmetric graph.
Since $\Gamma$ is $G$-vertex transitive
$(G, V \Gamma)$ is a transitive permutation group.
Since $A \Gamma \subseteq V \Gamma \times V \Gamma$
and $\Gamma$ is $G$-arc transitive, $A \Gamma$
is an orbital of $(G,V \Gamma)$.
It follows that $\Gamma = \orb_{A \Gamma}(G,V \Gamma)$.
\end{proof}
\end{proposition}

Relaxing the condition that $\Delta$ must be self-paired 
leads to a symmetric digraph.

\section{Imprimitive Symmetric Graphs}

Recall that if $\Delta$ is a subset of $\Omega$ 
and $g$ is an element of $G$, then $\Delta^g$ denotes
the image of $\Delta$ under the action of $g$.

\begin{definition}[Block of Imprimitivity]
If $(G,\Omega)$ is a transitive permutation group
then a subset $\Delta$ of $\Omega$
is said to be a \emph{block of imprimitivity} if for every $x \in G$
either $\Delta^x = \Delta$ or $\Delta^x \cap \Delta = \emptyset$
\end{definition}

\begin{definition}[$G$-invariant partition]
If $(G, \Omega)$ is a transitive permutation group,
then a partition $\cB$ of $\Omega$ is said to be \emph{$G$-invariant}
if for each $\Delta \in {\cB}$ and each $x \in G$
we have $\Delta^x \in {\cB}$. 
That is, $\cB$ admits $G$ as a group of permutations in a natural way.
\end{definition}

\begin{proposition}
If $(G,\Omega)$ is a transitive permutation group
and $\Delta$ is a block of imprimitivity
then $\cB = \{ \Delta^g : g \in G \}$ is a $G$-invariant
partition.
\begin{proof}
Since $(G,\Omega)$ is transitive,
every $\alpha \in \Omega$ is contained in $\Delta^g$
for some $g \in G$. Thus $\bigcup_{g \in G} \Delta^g = \Omega$.
If $\Delta^{g_1} \cap \Delta^{g_2} \neq \emptyset$
for some $g_1, g_2 \in G$
then $\Delta \cap \Delta^{g_2 g_1^{-1}} \neq \emptyset$.
Since $\Delta$ is a block of imprimitivity,
this implies that $\Delta = \Delta^{g_2 g_1^{-1}}$
and so $\Delta^{g_1} = \Delta^{g_2}$.
Thus $\cB$ is a partition of $\Omega$.
For any $\Delta^g \in \cB$ and any $x \in G$
$(\Delta^g)^x = \Delta^{gx} \in \cB$.
Thus $\cB$ is a $G$-invariant partition of $\Omega$.
\end{proof}
\end{proposition}

\begin{proposition}
If $(G,\Omega)$ is a transitive permutation group
and $\cB$ is a $G$-invariant partition of $\Omega$,
then each $B$ in $\cB$ is a block of imprimitivity.
\begin{proof}
Since $\cB$ is $G$-invariant, $B^g \in \cB$ for any
$g \in G$. Since $\cB$ is a partition, if $B^g \neq B$
then $B^g \cap B = \emptyset$.
The result follows.
\end{proof}
\end{proposition}

\begin{proposition}
The map $\pi : \Omega \to \cB$ which sends each point of $\Omega$
to the block of $\cB$ containing it
induces a surjective permutation homomorphism from $(G, \Omega)$
to $(G,\cB)$.
\begin{proof}
Need to show that $\pi(\alpha^x) = \pi(\alpha)^x$
for any $\alpha \in \Omega$ and any $x \in G$.
By the previous proposition $\pi(\alpha)$ is a block of imprimitivity.
Let $\Delta = \pi(\alpha)$. 
Since $\alpha \in \Delta$ we must have $\alpha^x \in \Delta^x$
That is $\pi(\alpha^x) = \Delta^x = \pi(\alpha)^x$.
\end{proof}
\end{proposition}

For any $B \in \cB$ we refer to the set $\pi^{-1}(B)$ as the
{\it fiber} of the homomorphism $\pi$ at $B$.
The fibers of a permutation homomorphism are blocks of imprimitivity.
Conversly each block of imprimitivity $\Delta$ 
gives rise to the $G$-invariant partition:
$\cB = \{ \Delta^x : x \in G \} $
and thus also to the homomorphism $(G,\Omega) \mapsto (G,\cB)$.

For every permutation group both $\{a\}$ and $\Omega$ are 
trivially blocks of imprimitivity. A permutation group is \emph{primitive}
if it admits no nontrivial blocks of imprimitivity.
%In category theoretic terms, the primitive representations
%of a group $G$ are those objects of $G$-Set having no proper quotients objects.
%Thus we see that the notions of transitivity and primitivity
%are category theoretic duals.

\begin{definition} [Imprimitive Symmetric Graph]
A symmetric graph $(G,\Gamma)$ is said to be {\it imprimitive} if
the induced permutation group, $(G, V \Gamma)$, is imprimitive.
\end{definition}

\begin{definition} [Quotient Graph]
Suppose that $(G,\Gamma)$ is an imprimitive symmetric graph
with $\cB$ a nontrivial $G$-invariant partition of the vertex set $V \Gamma$.
Define the {\it quotient graph} $\Gamma_\cB$ of $(G,\Gamma)$ with respect to $\cB$
to be the graph with vertex set $\cB$, and an arc $(B,C)$,
whenever there is some $\alpha \in B$ 
and some $\beta \in C$ 
such that  $(\alpha, \beta )$ is an arc of $\Gamma$.
\end{definition}

\begin{proposition}
$(G, \Gamma_\cB)$ is a symmetric graph.

\begin{proof}
We show that $G$ acts transitively on the arcs of $\Gamma_\cB$.
Let $(B,C)$ and $(D,E)$ be any two arcs of $\Gamma_\cB$.
By the definition of the quotient graph,
there must be some $\alpha \in \pi^{-1}(B)$ 
and $\beta \in \pi^{-1}(C)$
such that $(\alpha, \beta) \in A \Gamma$.
Simmilarly, there must be some $\gamma \in \pi^{-1}(D)$ 
and $\delta \in \pi^{-1}(E)$
such that $(\gamma, \delta) \in A \Gamma$,
By the symmetry of $(G, \Gamma)$, we can find an element $g \in G$
which carries  $(\alpha,\beta)$ to $(\gamma,\delta)$.
Since $\pi$ is a $G$-homomorphism, it follows that:
\begin{eqnarray*}
(B,C)^g & = & (\pi(\alpha),\pi(\beta))^g \\
        & = & (\pi(\alpha)^g,\pi(\beta)^g) \\
	& = & (\pi(\alpha^g),\pi(\beta^g)) \\
        & = & (\pi(\gamma),\pi(\delta)) \\
	& = & (D,E).
\end{eqnarray*}
\end{proof}
\end{proposition}

Suppose that some arc of $\Gamma$ 
has both its endpoints in the same fiber of $\cB$,
that is, there exists some $(\alpha,\beta) \in A \Gamma$
such that $\pi(\alpha) = \pi(\beta)$.
Since the fibers are blocks of imprimitivity,
it follows that $\pi(\alpha^g) = \pi(\beta^g)$ for each $g \in G$,
The arc transitivity of $\Gamma_\cB$ then implies 
every arc of $\Gamma$ has both its endpoints in the same fiber. 
In this case $\Gamma_\cB$ is the empty graph (no edges) 
with one vertex per connected component of $\Gamma$. 
This case is not very interesting,
to exclude it we say that the quotient is {\it nontrivial}
if it has valency at least one.
By the above discussion, $\Gamma_\cB$ is nontrivial
if and only if each of the fibers of $\cB$ is an
{\it independent set} of $\Gamma$.
For the remainder of this paper we shall always assume
that the quotient of a $G$-symmetric graph homomorphism
is non-trivial, even if we neglect to state this explicitely.

%Butterflies in bloom. See how the grow from a small seed into such majestic beings?

%% file: coset.tex
\chapter{Coset Graphs}

In this chapter we first show that for any group $G$, 
the transitive representations of $G$
together with $G$-homomorphisms between them
form a lattice isomorphic to a quotient of
the subgroup lattice of $G$.
%The primitive representations of $G$ are the co-atoms
%of this lattice.
Next we describe a construction of Sabidussi's
for vertex transitive graphs, and give a 
group theoretic characterization
of symmetric graphs.
Finally we consider quotients 
of symmetric graphs from a group theoretic perspective. 

%We describe three well-known methods of constructing ``extensions''
%of symmetric graphs -- the direct sum, the lexicographic product,
%and the Biggs covering graph construction.

\section{Transitive Permutation Groups}

% Coset Representation
\begin{definition}[Core of a subgroup]
If $H$ is a subgroup of $G$, then the \emph{core} of $H$ in $G$ is:
$$\Core_G(H) = \bigcap_{x \in G} x^{-1} H x $$
\end{definition}

\begin{proposition} If $H$ is a subgroup of $G$, 
then the core of $H$ in $G$ is a normal subgroup of $G$.
\begin{proof}
Since $gx$ runs over all the elements of $G$ as $x$ does, we have:
\begin{eqnarray*}
g^{-1}\Core_G(H)g & = & \bigcap_{x \in G} (gx) H (gx)^{-1} \\
                  & = & \bigcap_{x \in G} x H x^{-1} 
\end{eqnarray*}               
for any $g \in G$. 
\end{proof}
\end{proposition}

\begin{definition} [Coset Representation]
For $H$ a subgroup of $G$,
let $\Cos_G(H)$ denote the right cosets of $H$ in $G$.
We may define an action of $G$ on $\Cos_G(H)$ by:
$$ (Ha)^x = Hax $$
This is indeed an action, since:
\begin{eqnarray*}
{(Ha)^x}^y & = & (Hax)^y \\
           & = & (Haxy) \\
	   & = & (Ha)^{xy}
\end{eqnarray*}
for any ``$Ha$'' a coset of $H$ in $G$, and any $x,y \in G$.
We call a permutation group of the form $(G,\cos_G(H))$ a
{\it coset representation} of $G$.
\end{definition}

\begin{example}[Right Regular Representation]
For any group $G$, the \emph{right regular representation}
of $G$ is the permutation group $(G,G)$ with the action given by:
$$ g^{h} = gh. $$
This is a special case of a coset representation 
where $H$ is the trivial group $\{ 1 \}$.
\end{example}

\begin{proposition}
For any pair of groups $H \leq G$, the action of $G$
on the cosets of $H$ is transitive with kernel $\Core_G(H)$.
\begin{proof}
Let $Ha_1$ and $Ha_2$ be any two cosets of $H$ in $G$,
then:
$$(Ha_1)^{a_1^{-1}a_2} = Ha_2$$ 
so the action is transitive.

The stabilizer of the point ``$H$'' is $H$.
By Lemma 1, Chapter 1, 
the stabilizer of the point ``$Ha$'' is $a^{-1}Ha$.
If $g \in G$ is such that $g$ stabilizes every coset of $H$ in $G$,
then we must have $g \in a^{-1}Ha$ for each $a \in G$.
That is $g \in \Core_G(H)$.

\end{proof}
\end{proposition}

The next proposition tells us that every transitive representation
of a group $G$ is permutation equivalent to some coset representation.

\begin{proposition}
If $(G,\Omega)$ is a transitive permutation group
and $\alpha$ is some point of $\Omega$, then
$(G, \Omega)$ is permutation equivalent to $(G, \cos_G(G_\alpha))$

\begin{proof}
For each $\beta \in G$ Consider:
$$ S_\beta = \{ x \in G : \alpha^x = \beta  \} $$
If $\beta = \alpha^g$ then it's not hard to see
that $S_{\beta}$ is a right coset of $G_{\alpha}$:
$$ S_{\beta} = G_{\alpha} g $$
Define $\eta : \Omega \to \cos_G(G_a)$ by $\eta(\beta) = S_\beta$
and take $\varphi : G \to G$ to be the identity.
We have, for any $\beta \in \Omega$ and any $x \in G$,
if $\beta = \alpha^g$ then:
\begin{eqnarray*}
\eta(\beta^x) & = & S_{\beta^x} \\
              & = & S_{\alpha^{gx}} \\
	      & = & G_\alpha gx \\
	      & = & G_\alpha g^{\varphi(x)} \\
	      & = & \eta(\alpha^g)^{\varphi(x)} \\
	      & = & \eta(\beta)^{\varphi(x)}
\end{eqnarray*}
Thus $\eta$ and $\varphi$ establish the desired permutation equivalence.
\end{proof}
\end{proposition}
	
\section{Imprimitive Permutation Groups}
	
\begin{definition}
If $(G,\Omega)$ is a transitive permutation group and $\Delta$
is a subset of $\Omega$ then the {\it setwise stabilizer} of $\Delta$ is:
$$G_{\Delta} = \{ x \in G : \Delta^x = \Delta \}.$$
\end{definition}

\begin{lemma}
If $(G,\Omega)$ is a transitive permutation group and $\Delta$
is a block of imprimitivity, then
$G_{\alpha} \leq G_{\Delta} \leq G$.
for any $\alpha \in \Delta$.
\begin{proof}
Since $\alpha \in \Delta$, 
for any $g \in G$ we have $\alpha^g \in \Delta^g$.
Since $\Delta$ is a block of imprimitivity,
if $\alpha^g = \alpha$, then $\Delta^g \cap \Delta \neq \emptyset$
so $\Delta^g = \Delta$. That is $G_\alpha \leq G_\Delta$.
\end{proof}
\end{lemma}

\begin{proposition}
If $(G,\Omega)$ is a transitive permutation group and $\Delta$
is a block of imprimitivity, then $(G_{\Delta},\Delta)$
is a transitive permutation group.
\begin{proof}
Suppose that $\alpha, \beta \in \Delta$.
Since $(G,\Omega)$ is transitive, there is some $g \in G$
such that $\alpha^g = \beta$. 
Since $\alpha^g \in \Delta^g$ and $\beta \in \Delta$,
we must have $\Delta^g \cap \Delta \neq \emptyset$.
But $\Delta$ is a block of imprimitivity,
so this implies that $\Delta^g = \Delta$. 
That is $g \in G_{\Delta}$,
so $G_{\Delta}$ acts transitively on $\Delta$.
\end{proof}
\end{proposition}

\begin{proposition}
Suppose that $(G,\Omega)$ is a permutation group and $\alpha \in \Omega$.
For any subgroup $H$ such that $G_\alpha \leq H \leq G$, 
the set $\alpha^H = \{\alpha^h : h \in H \}$ is a block of imprimitivity.
\begin{proof}
Let $\Delta = \alpha^H$.
If $\Delta^g \cap \Delta \neq \emptyset $ for some $g \in G$
then there must be some $ x,y \in H $ such that
$\alpha^{xg} = \alpha^y$.
It follows that $ g \in x^{-1} G_\alpha y \leq H $
and $ \Delta^g = \Delta $.
That is $ \Delta^g \cap \Delta \neq \emptyset \Rightarrow \Delta^g = \Delta$,
so $\alpha^H$ is a block of imprimitivity. 
\end{proof}
\end{proposition}

\begin{lemma}
If $(G, \Omega)$ is a transitive permutation group
and $\alpha$ is some point of $\Omega$,
then $G_{\alpha^H} = H$
for any $H$ such that $G_\alpha \leq H \leq G$.
\begin{proof}
Let $\Delta = \alpha^H$.
Clearly $H \leq G_{\Delta}$.
To show that $H = G_\Delta$ it suffices to show that $G_\alpha$
and $G_\Delta$ have the same index in $H$.
Since $\alpha^H$ is an orbit of $H$ on $\Omega$,
the permutation group $(H,\alpha^H)$ is transitive.
Since $G_\alpha \leq H$, 
the stabilizer of the point $\alpha$ in $(H,\alpha^H)$
is $H \cap G_\alpha = G_\alpha$.
Thus by the orbit stabilizer theorem we have
$\Delta = [H : G_\alpha]$.

On the other hand, by proposition 11 and 12,
$\Delta$ is a block of imprimitivity
and $(G_{\Delta}, \Delta)$ is a transitive permutation group.
By lemma 3 we have $G_\alpha \leq G_\Delta$ 
so the stabilizer of the point $\alpha$ in $(G_\Delta, \Delta)$
is $G_\Delta \cap G_\alpha = G_\alpha$.
So again by the orbit stabilizer lemma we have
$|\Delta| = [G_{\Delta}: G_{\alpha}]$.
The result follows.

\end{proof}
\end{lemma}

\begin{lemma}
If $(G,\Omega)$ is a transitive permutation group,
$\Delta$ is a block of imprimitivity,
and $\alpha$ is any point of $\Delta$,
then $\alpha^{G_{\Delta}} = \Delta$.
\begin{proof}
Clearly $\alpha^{G_{\Delta}} \leq \Delta$.
Let $H = G_{\Delta}$.
We have $G_{\alpha} \leq H \leq G$
and $(H,\alpha^H)$ is a transitive permutation group
equivalent to $(H,\Cos_H(G_\alpha))$.
Since $(G_{\Delta}, \Delta)$ is also a transitive permutation group
equivalent to $(H, \Cos_H(G_\alpha))$
we must have $|Delta| = |alpha^H|$.
Since $\alpha^H \leq \Delta$, this implies that $\alpha^H = \Delta$.

\end{proof}

\end{lemma}

Let $(\cS, \leq)$ be the set of subgroups of $G$ containing $G_\alpha$,
partially ordered by the subgroup relation.
Let $(\cP, \subseteq)$ be the set of blocks of imprimitivity containing $\alpha$,
partially ordered by the subset relation.

\begin{proposition}
$(\cS, \leq)$ is order isomorphic to $(\cP, \subseteq)$

\begin{proof}
Let $\Phi : \cS \to \cP$ be given by
$\Phi(\Delta) = G_{\Delta}$
and let $\Psi : \cP \to \cS$ be given by
$\Psi(H) = \alpha^H$
Making use of lemma ?? and ??,
for any $H \in \cS$ we have:
$$ \Phi \circ \Psi(H) = \Phi(\alpha^H) = G_{\alpha^H} = H $$
and for any $\Delta \in \cP$ we have:
$$ \Psi \circ \Phi(\Delta) = \Psi (G_{\Delta}) = \alpha^{G_{\Delta}} = \Delta $$
so $\Phi$ and $\Psi$ are inverses.
To see that $\Phi$ is order preserving note that:
$$ G_{\Delta_1} \leq G_{\Delta_2} 
\Longleftrightarrow \alpha^{G_{\Delta_1}} \leq \alpha^{G_{\Delta_2}}
\Longleftrightarrow \Delta_1 \subseteq \Delta_2  $$
\end{proof}
\end{proposition}

\begin{proposition}
Let $(G,\Omega_1)$ and $(G,\Omega_2)$ be two transitive representations
of $G$. For any $\alpha \in \Omega_1$, $\beta \in \Omega_2$,
the permutation groups $(G,\Omega_1)$ and $(G,\Omega_2)$
are permutation equivalent if and only if there is an automorphism of $G$
carrying $G_\alpha$ to $G_\beta$.
\begin{proof}
Suppose that $(G, \Omega_1)$ and $(G, \Omega_2)$ are permutation equivalent.
Let $\varphi: G \to G$ and $\eta : \Omega_1 \to \Omega_2$ be the maps
establishing the equivalence. 
Let $\gamma = \eta(\alpha)$.
Clearly $G_\gamma = \varphi(G_\alpha)$.
By transitivity there must be some $g \in G$ such that $\beta = \gamma^g$.
Let $\vartheta_g : G \to G$ be the map: $x \mapsto g^{-1} x g$.
We have:
\begin{eqnarray*} 
G_\beta & = & g^{-1}G_\gamma g \\
        & = & \vartheta(G_\gamma) \\
	& = & \vartheta \circ \varphi (G_\alpha)
\end{eqnarray*}
So $\vartheta \circ \varphi$ carries $G_\alpha$ to $G_\beta$.

For the other direction, let $\psi$ be the map carrying $G_\alpha$
to $G_\beta$. By proposition 1,
$(G, \Omega_1)$ is permutation equivalent to $(G, \cos(G_\alpha))$
and $(G, \Omega_2)$ is permutation equivalent to $(G, \cos(G_\beta))$.
Take $\varphi = \psi$ and let 
$\eta: \cos(G_\alpha) \to \eta(G_\beta)$ be given by
$(G_\alpha) g \mapsto (G_\beta) \psi(g)$.
For any $x \in G$ we have:
\begin{eqnarray*}
\eta({(G_\alpha g)}^x) & = & \eta((G_\alpha) gx) \\
                     & = & (G_\beta) \psi(gx)  \\
		     & = & {(G_\beta) \psi(g)}^{\psi(x)} \\
		     & = & \eta((G_\alpha) g)^{\varphi(x)}
\end{eqnarray*}
thus $\varphi$ and $\eta$ establish a permutation equivalence between
$(G, \cos(G_\alpha)$ and $(G, \cos(G_\beta)$ and hence between
$(G, \Omega_1)$ and $(G, \Omega_2)$.
\end{proof}
\end{proposition}

Define an equivalence relation on the set of subgroups of $G$
as follows: $H_1 \sim H_2$ if and only if there is an automorphism of $G$
carrying $H_1$ to $H_2$. 
Let $\cZ$ denote the equivalence classes of this relation.
Define a partial order on $\cZ$
by $[H_1] \leq [H_2]$ if and only if there is some $\tilde{H_1} \in [H_1]$ and some
$\tilde{H_2} \in [H_2]$ such that $\tilde{H_1} \leq \tilde{H_2}$. 

The previous result tells us that, upto permutation equivalence,
the essential information about the transitive representations of $G$ 
is contained in $(\cZ, \leq)$. In particular all the transitive representations
of $G$ are to be found as quotients of the right regular representation.

\section{Sabidussi's Construction}

% Better make sure to discuss the difference between
% isomorphic as graphs and isomorphic as $G$-graphs.

We saw in the last section that
every permutation group is equivalent to
one of the form $(G,\cos_G(H))$.
Thus, upto isomorphism, transitive permutation groups are 
uniquely determined by pairs of groups $(G,H)$ 
with $H$ a subgroup of $G$.
In this section and the next
we shall see that upto isomorphism,
a symmetric graph is uniquely
determined by a triple $(G,H,a)$ 
where $H$ is a subgroup of $G$ 
and $a \in G \backslash H$ is an involution.

Given a triple $(G,H,a)$, the idea is to construct
a graph whose vertices are the cosets of $H$ in $G$.
We call such a graph a {\it coset graph}.
The idea of a coset graph originally goes back to Sabidussi
who was studying vertex transitive graphs.
Sabidussi's construction was essentially 
a generalization of the {\it Cayley graph} construction.

\begin{definition} (Cayley Graph)
Given a group $G$ and a subset $D \subseteq G$ the {\it Cayley graph},
$\Cay(G,D)$ is the directed graph with vertices the elements of $G$
and arc set $\{(x,y) : xy^{-1} \in D \}$
\end{definition}

A Cayley Graph will contain no loops,
provided that $1 \not\in D$.
A Cayley Graph is simple if and only if 
the set $D$ is closed under taking inverses,
and connected if and only if $D$ is a generating set for $G$.
The vertices are the elements of the group $G$,
and the action of $G$ on the vertices is permutation equivalent 
to the right regular representation of $G$.

Not every vertex transitive graph is a Cayley Graph.
Sabidussi's idea was to generalize Cayley's construction by considering the 
action of $G$ on the cosets of some non-trivial subgroup $H$.

\begin{definition} (Sabidussi Graph)
Given a group $G$, a subgroup $H$, and a set $D \subseteq G$, 
the Sabidussi graph, $\Sab(G,H,D)$, 
is the directed graph with vertex set $\cos_G(H)$ 
and arc set $\{(Hx,Hy) : xy^{-1} \in D\}$
\end{definition}

A Sabidussi graph will contain no loops
provided that $D \cap H = \emptyset$.
A Sabidussi Graph is simple if and only if
$D$ is closed under taking inverses
and connected if and only if $D \cup H$ generates $G$.
The vertices of the Sabidussi graph are  
the cosets of the subgroup $H$ in $G$,
and the action of $G$ on the vertices is $(G,cos(H))$.

\begin{proposition}[Sabidussi]
Sabidussi Graphs are vertex transitive.
Every vertex transitive graph is isomorphic to a Sabidussi Graph.
\begin{proof}
Let $\Gamma = \Sab(G,H,D)$ be any sabidussi graph
(simple and without loops).
The group $G$ has a natural transitive action on the vertices of $\Gamma$
given by $Hx^g = Hxg$.
To prove that $\Gamma$ is $G$-vertex transitive,
we must check that this action preserves the adjacency structure of $\Gamma$.

Suppose that $(Hx,Hy)$ is an arc of $\Gamma$,
so $xy^{-1} \in D$.
Since $xg(yg)^{-1} = xgg^{-1}y^{-1} = xy^{-1} \in D$
for any $g$, 
it follows immediately that $(Hx,Hy)^g = (Hxg,Hyg)$
is also an arc of $\Gamma$.
This proves the first part.

Now, let $\Gamma$ be any $G$-vertex transitive graph.
It follows that $(G, V \Gamma)$ is permutation
equivalent to $(G, \cos_G(H))$ for some $H < G$.
Let $\eta : \Omega \to \cos_G(H)$ be the map establishing
the permutation equivalence, and let $\mu = \eta^{-1}$

Let $N$ be the neighbourhood of the vertex $\mu(H)$, and let 
$$D = (\bigcup_{\alpha \in N} \eta(\alpha)) \diagdown H$$
We will show that $\Gamma \cong \Sab(G,H,\eta(D))$.

The map $\eta : V \Gamma \to \cos_G(H)$
establishes a permutation isomorphism between the
vertices of $\Gamma$ and the vertices of $\Sab(G,H,D)$.
Let $(\alpha, \beta)$ be any arc of $\Gamma$.
By $G$-vertex transitivity, we can find some $g \in G$
such that $\alpha^g = \mu(H)$ and $\beta^g \in E$.
Suppose that $\eta(\alpha) = Hx$ and $\eta(\beta) = Hy$.
We must check that $xy^{-1} \in D$

Since $\eta$ is a permutation isomorphism,
it follows that $\eta(\alpha^g) = Hxg$
and $\eta(\beta^g) = Hyg$.
Since $\beta^g \in E$ we must have 
$xg(yg)^{-1} = xgg^{-1}y = xy^{-1}\in D$.
The result follows.

\end{proof}

\end{proposition}

\section{A Group Theoretic Characterization of Symmetric Graphs}

Since symmetric graphs are vertex transitive,
every symmetric graph is a Sabidussi graph.
However, not every vertex transitive graph is symmetric,
so we expect there to be some extra conditions on
the subset $D$ in the symmetric case.

Recall from Chapter 1 that every symmetric graph 
is an orbital graph.
In the first section of this chapter, we saw that if $H < G$
is the stabilizer of some point in $\Omega$, then $(G,\Omega)$
is permutation equivalent to $(G,\cos_G(H))$.
We shall see in this section, the orbitals of $(G,\Omega)$
actually correspond to {\it double cosets} of $H$ in $G$

\begin{definition}
If $H$ is a subgroup of $G$, then a {\it double coset} of $H$ in $G$
is a subset of $G$ of the form
$HxH = \{ h_1 x h_2 : h_1, h_2 \in H \}$
for some $x \in G$.
\end{definition}

\begin{lemma} 
Each double coset of $H$ in $G$ is a union of right cosets of $H$.
The double cosets of $H$ in $G$ form a partition of $G$.
\begin{proof}
The first part is obvious since $HxH = \bigcup_{h \in H} Hxh$
for any $x \in G$.
For each element $g \in G$ it is clear that $g \in HgH$
so $\bigcup_{x \in G} HxH = G$.
If $HxH \bigcap HyH \neq \emptyset$ 
then $h_1 x h_2 = h_3 y h_4 $ for some $h_1,h_2,h_3,h_4 \in H$.
It follows that $y = h_1^{-1} h_3 x h_4 h_2^{-1}$ 
and $HyH = H h_1^{-1} h_3 x h_4 h_2^{-1} H = HxH$.
This shows that any two double cosets are equal or disjoint,
so the double cosets partition $G$ as claimed.
\end{proof}
\end{lemma}

\begin{proposition}[Lorimer]
For any permutationg group $(G,\Omega)$ with point stabilizer $H$, 
the orbitals are in natural bijection with the double cosets of $H$ in $G$.
The self-paired orbitals correspond to those 
double cosets which contain an involution.
\begin{proof}
From the first section in this chapter we know that $(G,\Omega)$
is permutation equivalent to $(G,\cos_G(H))$.
Recall from chapter 1 that the orbitals of a permutation group
are in bijection with the orbits of the point stabilizer.
We must show that orbit of $(H,\cos_G(H))$
corresponds to a double coset of $H$ in $G$.

Let $\kappa$ be the map $Hx \mapsto HxH$.
Suppose that $Hx$ and $Hy$ are in the same orbit
of $(H,\cos_G(H))$.
Then there is some $h \in H$ such that $Hxh = Hy$.
It follows that $y \in Hxh \subset HxH$.
That is $HyH = HxH$ and so $\kappa(Hx) = \kappa(Hy)$.

Conversly, suppose that $\kappa(Hx) = \kappa(Hy)$
for some pair of cosets $Hx$ and $Hy$.
Then $x \in HyH$ and so $x = h_1 y h_2$ for some 
$h_1, h_2 \in H$.
It follows that $Hx = Hyh_2$,
that is $Hx$ and $Hy$ lie in the same orbit of 
$(H,\cos_G(H))$.
Since $\kappa$ is surjective,
this proves the first part. 

Suppose that $a \in G$ is such that $a^2 = 1$ and $a \not\in H$.
Since $(H, Ha)^a = (Ha, Ha^2) = (Ha, H)$, it follows that 
the orbital $\Delta$ of $(G,\cos_G(H))$ 
containing $(H,Ha)$ is self-paired.
Clearly $a \in HaH = \kappa(Ha)$,
so the double coset associated with $\Delta$ contains an involution.

Conversely, if the double coset $HxH$ contains 
an involution $a$, then $HxH = HaH = \kappa(Ha)$.
Thus $HxH$ is associated with the orbital containing
$(H,Ha)$ which is self-paired.
This proves the second part.
\end{proof}
\end{proposition}

\begin{theorem} [Lorimer] For any pair of groups $H \leq G$ and any $a \in G$
such that $a \not\in H$ and $a^2 = 1$, the Sabidussi graph
$\Sab(G,H,HaH)$ is $G$-symmetric. Furthermore every symmetric graph
is of this form for some $G$, $H$ and $a$.
\begin{proof} 
Let $\Gamma = \Sab(G,H,HaH)$.
By proposition ?? $\Gamma$ is $G$-vertex transitive.
Since $HaH$ is an orbit of the action of $H$
on $\Cos_G(H)$, the graph $\Gamma$ is also $G$-locally transitive.

Suppose that $\Gamma$ is a $G$-symmetric graph.
Let $\alpha$ be any vertex, and let $H = G_\alpha$.
By proposition ??, the permutation group $(G, V \Gamma)$ 
is equivalent to $(G, \Cos_G(H))$.
Let $\eta : V \Gamma \to \Cos_G(H))$ be any map
inducing this equivalence.
Let $\beta$ be any neighbour of $\alpha$ in $\Gamma$
and let $a \in G$ be such that $(\alpha,\beta)^a = (\beta,\alpha)$.
Clearly $a$ is an involution and $\eta(\beta) = Ha$.

Let $\Gamma' = \Sab(G,H,HaH)$.
The map $\eta$ induces a permutation equivalence
between the vertices of $\Gamma$ and the vertices of $\Gamma'$.
To show that $\Gamma$ is isomorphic to $\Gamma'$
we must show that $\eta$ preserves adjacency.

Let $(\omega, \delta)$ be any arc of $\Gamma$.
Let $g \in G$ be such that $\omega^g = \alpha$.
Then $\delta^g$ is some neighbour of $\alpha$.
Thus there exists $h \in H$ such that $\delta^{gh} = \beta$.
Therefore $\eta(\delta^{gh}) = \eta(\beta) = Ha$.
Since $\eta$ is a permutation equivalence,
$\eta(\delta^{gh}) = \eta(\delta)^{gh}$,
so $\eta(\delta) = Ha^{{(gh)}^{-1}} = Hah^{-1}g^{-1}$.
Simmilarly, $\eta(\omega^g) = \eta(\omega)^g$
so $\eta(\omega) = H^{g^{-1}} = Hg^{-1}$.
Now $ah^{-1}g^{-1}g = ah^{-1} \in HaH$
so $(\eta(\omega),\eta(\delta))$ is an arc of $\Gamma'$.
The result follows.

\end{proof}
\end{theorem}

%\begin{proposition}
%$\Gamma = \Sab(G,H,HaH)$ is connected if and only if together
%$H$ and $a$ generalte $G$
%\end{proposition}

\begin{lemma}
$|HxH| = |H||H|/|x^{-1}Hx \bigcap H|$ Fix me
\begin{proof}
The mapping $HxH \to x^{-1}HxH$ given by $h_1 x h_2 \mapsto x^{-1} h_1 x h_2$
      is bijective, so $|HxH| = |x^{-1}HxH|$. 
      But for the rule for the product of two subgroups: 
      $|x^{-1}HxH| = |x^{-1}Hx||H|/|x^{-1}Hx \bigcap H|$,
      of course, $|x^{-1}Hx| = |H|$. and the result follows.
\end{proof}
\end{lemma}

\begin{proposition}
If $\Gamma$ is a $G$-symmetric graph isomorphic to the Sabidussi
graph $Sab(G,H,HaH)$ the the stabilizer of an arc of $\Gamma$
is isomorphic to $a^{-1}Ha \cap H$ and the valency of $\Gamma$
is $ |H|/|x^{-1}Hx \bigcap H| $.
\begin{proof}
Consider the arc $(H,Ha)$. The stabilizer of this arc is
the intersection between the stabilizer of the vertex ``$H$''
and the stabilizer of the vertex ``$Ha$''.
But the stabilizer of the vertex ``$H$'' is just $H$
and the stabilizer of the vertex ``$Ha$'' is $a^{-1}Ha$.
So the stabilizer of the arc is $a^{-1}Ha \cap H$.

Consider the action of $H$ on the neighbours 
of the vertex ``$H$''.
Since $\Gamma$ is symmetric, this action is transitive.
If ``$Ha$'' is any neighbour of ``$H$'',
then the action of $H$ on the neighbours of ``$H$''
is permutation equivalent to the action of $H$
on the cosets of the stabilizer of the arc $(H,Ha)$.
Since the stabilizer of the arc $(H,Ha)$ is $a^{-1}Ha \cap H$
it follows that the valency of $\Gamma$ is
$[H : a^{-1}Ha \cap H] = |H|/|x^{-1}Hx \bigcap H|$

\end{proof}

\end{proposition}

In some sense the valancy of a symmetric graph is
measuring the extent to which a subgroup fails to be normal.

\section{Quotient Graphs}

In the first section of this Chapter,
we determined, for any two transitive permutation groups
$(G, \Omega_1)$ and $(G,\Omega_2)$, 
the conditions under which there exists a
$G$-homomorphism from $(G,\Omega_1)$ to $(G,\Omega_2)$.

We saw that if $(G,\Omega_1)$
is permutation equivalent to $(G,\cos_G(H))$ 
for some $H < G$,
and $(G,\Omega_2)$ is permutation equivalent
to $(G,\cos_G(K))$ for some $K < G$,
then there is a $G$-homomorphism from $(G,\Omega_1)$
to $(G,\Omega_2)$ if and only if $H < K < G$.

In this section we give analogous conditions
under which there exists a $G$-homomorphism
from a $G$-symmetric graph $\Gamma$
to a $G$-symmetric graph $\Sigma$.

\begin{theorem} [Lorimer]
If $\Gamma$ is a $G$-symmetric graph isomorphic to $\Sab(G,H,HaH)$
and $\Sigma$ is a quotient of $\Gamma$,
then $\Sigma$ is isomorphic to $\Sab(G,K,KaK)$ for some $H < K < G$.
\begin{proof}
Let $\eta : V \Gamma \to \Cos_G(H)$ be the map inducing
the isomorphism between $\Gamma$ and $\Sab(G,H,HaH)$.
Let $\pi : V \Gamma \to V \Sigma$ be the map
inducing the homomorphism between $\Gamma$ and $\Sigma$.
Let $\alpha = \eta^{-1}(H)$,
let $B = \pi(\alpha)$ and let $\Delta = \pi^{-1}(B)$.

By proposition ?? $\Delta$ is a block of imprimitivity.
If $K = G_{\Delta}$ then the permutation group $(G, V \Sigma)$
is equivalent to $(G, \cos_G(K))$.
Let $\mu : V \Sigma \to \cos_G(K)$ be any map
inducing the equivalence.
If $\Sigma' = \Sab(G,K,KaK)$ then $\mu$ induces
a permutation equivalence between the vertices of $\Sigma$
and the vertices of $\Sigma'$.
We must check that it preserves adjacency. Blah blah blah...

\end{proof}
\end{theorem}

If $a \in K$ then the valency of the quotient is one.

%Butterflies in bloom. See how the grow from a small seed into such majestic beings?

%% file: geometric.tex
\chapter{The Extension Problem}

We would like to understand
how a symmetric graph $(G, \Gamma)$ can be ``unfolded''
into a larger imprimitive symmetric graph $(\widetilde{G}, \widetilde{\Gamma})$
admitting the original graph as a quotient.
In particular we'd like to understand how pairs of graphs
$(\Gamma, \widetilde{\Gamma})$ where $\widetilde{\Gamma}$
is an ``extension'' of $\Gamma$ are related combinatorially.
That is, we'd like to be able to describe the the structure
of the graph $\widetilde{\Gamma}$
in terms of the structure of the graph $\Gamma$.

The quotient of a symmetric graph contains considerably less information
than the original graph.
Gardiner and Praeger \cite{GPL1993} observed that
some of the information that is lost
may be recovered from the induced bipartite graph 
between adjacent blocks of the partition,
and from a combinatorial design induced on the blocks themselves.

In this chapter we describe Gardiner and Praeger's observations
as well as some of the questions is raises.

\section{Induced Bipartite graph}

% Introductory remarks

Let $\Gamma$ be an imprimitive $G$-symmetric graph,
with $\cB$ a non-trivial $G$-invariant partition of the vertices.
As always, we assume that the quotient $\Gamma_\cB$ has valency at least one,
so the blocks of $\cB$ are independent sets.
For any arc $(B,C)$ of the quotient,
the subgraph of $\Gamma$ induced by $B \cup C$ must be bipartite --
possibly containing some isolated vertices.
If we restrict ourselves to the subgraph induced by
$(\Gamma(C) \cap B) \cup (\Gamma(B) \cap C)$
then we obtain a bipartite graph with no isolated vertices.
We call this graph the {\it induced bipartite graph} of $(B,C)$,
and denote it by $\Gamma[B,C]$.

\begin{proposition}
The induced bipartite graph $\Gamma[B,C]$ is $G_{B \cup C}$ symmetric.
\begin{proof}
Let $(\alpha,\beta)$ and $(\gamma,\delta)$ be any two arcs of $\Gamma[B,C]$.
Without loss of generality we may assume that $\alpha \in B$ and $\beta \in C$.
By the $G$-symmetry of $\Gamma$ we can find $g \in G$ such that
$(\alpha,\beta)^g = (\gamma,\delta)$.
Since $\Gamma[B,C]$ is bipartite, either $\gamma \in B$ and $\delta \in C$ 
in which case $B^g = B$ and $C^g = C$,
or $\gamma \in C$, $\delta \in B$, in which case $B^g = C$ and $C^g = B$.
Either way, $g \in G_{B \cup C}$ and the result follows.
\end{proof}
\end{proposition}

\begin{proposition}
For any two arcs in the quotient, the induced bipartite graphs are isomorphic.
\begin{proof}
Let $(B,C)$ and $(D,E)$ be any two arcs of $\Gamma_\cB$.
Since $\Gamma_\cB$ is $G$-symmetric we can find some $g \in G$
such that $(B,C)^g = (D,E)$. Clearly
$$(\Gamma(C) \cap B)^g = (\Gamma(C^g) \cap B^g) = (\Gamma(E) \cap D),$$
and
$$(\Gamma(B) \cap C)^g = (\Gamma(B^g) \cap C^g) = (\Gamma(D) \cap E).$$
So $g$ induces a bijection
$$(\Gamma(C) \cap B) \cup (\Gamma(B) \cap C)
\to (\Gamma(E) \cap D) \cup (\Gamma(D) \cap E).$$
Since $g$ is an automorphism, adjacency is preserved and we have an isomorphism
from $\Gamma[B,C]$ to $\Gamma[D,E]$.
\end{proof}
\end{proposition}

%Note that if $(B,C)$ and $(D,E)$ are any two arcs of $\Gamma_\cB$
%then $G_{B \cup C}$ is isomorphic to $G_{D \cup E}$
%though they are not actually {\it the same} group.
%They are different isomorphic subgroups of $G$.
%Let $\cA$ denote the isomorphism class of groups
%to which $G_{B \cup C}$ and $G_{D \cup E}$ both belong.
%The two grphs $\Gamma[B,C]$ and $\Gamma[D,E]$
%are $\cA$-equivalent, but not $\cA$-isomorphic.

When passing from a symmetric graph to its quotient
we preserve the adjacency structure of the blocks,
and discard the exact details of which vertices are connected to which others.
Some of the information that has been lost is recoverable from the induced bipartite graph.
But the induced bipartite graph reveals only the {\it local} connectivity.
To reconstruct the original graph from its quotient we need to know
how these induced bipartite graphs are ``glued together''.
Gardiner and Praeger observed that some of this global information about
the way the induced bipartite graphs fit together is captured in a combinatorial design
induced on each of the blocks.
Before we can describe this design we need a little background.

\section{Combinatorial Designs}

% Need to talk about connection between design theory and geometry.

%Connection with group theory:
%Many interesting groups, including a number of
%the sporadic simple groups arise as automorphism
%groups of $t$-designs.

\begin{definition}
An {\it incidence structure} is a triple $(P,B,I)$
where $P$, and $B$ are sets, usually refered to as
{\it points} and {\it blocks} respectively,
and $I$ is an {\it incidence relation} $I \subseteq P \times B$.
\end{definition}

It is often convenient to visualize an incidence structure
as a bipartite graph. 
Take the points as vertices of one side of
the bipartition and the blocks as vertices of the other.
Draw an edge between the vertex corresponding to a point $p$
and the vertex corresponding to a block $b$ if and only if
$(p,b) \in I$.

We will be interested in incidence structures
satisfying strong regularity and symmetry conditions.

\begin{definition}
A $(v,k,\lambda)$-design is an incidence structure such that:
\begin{enumerate}
\item There are $v$ points in total
\item Each block is incident with exactly $k$ points
\item Each point is incident with exactly $\lambda$ blocks
\end{enumerate}
\end{definition}

When visualizing an incidence structure as a bipartite graph,
the extra regularity conditions of a design correspond 
to the condition that 
any two vertices in the same bipartite half of the graph
must have the same valency.

Note that our definition here of a ``design'' 
corresponds to what is usually refered to as a $1$-design 
or ``tactical configuration'' in the literature.
There is a more general definition of a $t$-design
of which the definition of a $1$-design is a special case.
We will not be needing this more general definition.

% Connection with geometry?

%As an example, the points and lines of the projective plane 
%over a field of $q$ elements, for $q$ a prime power, 
%form a $(q^2 + q + 1,q + 1,q + 1)$-design. 
%The points and lines of the affine plane over a field with $q$ elements 
%form a $(q^2,q,q)$ design. 
%In a very loose sense designs, particularly $2$-designs,
%may be viewed as generalized finite geometries.

The incident point-block pairs of a design are usually refered to as {\it flags}.
For a design $\cD$ we will use the notation
$P_\cD$, $B_\cD$ and $F_\cD$ to denote the points,
blocks and flags of $\cD$ respectively.

Define the {\it trace} of a block $b$ to be the set
$T(b) = \{ p \in P : (p,b) \in I \}$.
Simmilarly define the {\it trace} of a point $p$
to be the set $T(p) = \{b \in B : (p,b) \in I \}$.

\begin{proposition}
In a design, the number of blocks with the same trace is a constant 
that is independent of the choice of block
\end{proposition}

We denote this constant by $m$ and call it the 
{\it multiplicity} of the design.
If distinct blocks have distinct traces,
then we may identify the blocks with their traces,
and take $B$ to be a subset of the power set of $P$.
Otherwise we say that the design contains {\it repeated blocks}.

\begin{proposition}
If a design contains repeated blocks,
then the incidence structure obtained by identifying
blocks with the same trace is again a design
\end{proposition}

\begin{definition}
if $\cD = (P_\cD, B_\cD, I)$ is a $1$-design then the {\it dual} 
of $\cD$ is the design $ \cD^{*} = (B_\cD, P_\cD, I^{*})$
where $(b,p) \in I^{*}$ if and only if $(p,b) \in I$.
\end{definition}

\begin{proposition}
The dual of a $1$-design is a $1$-design.
\end{proposition}

\section{Flag Transitive Designs}

An automorphism of a design $\cD$ is a pair $(\eta, \mu)$ of bijections
$\eta : P_\cD \to P_\cD$ and $\mu : B_\cD \to B_\cD$
with the property that $(p,b)$ is a flag of $\cD$
if and only if $(\eta(p), \mu(b))$ is a flag.
The automorphisms of a design form a group.

As with sets and graphs, we say that a group $G$ acts on a design $\cD$
as a group of automorphisms if there is some homomorphism $\rho : G \to \Aut(\cD)$.
We don't require this homomorphism to be injective, that is we don't
require that $G$ acts {\it faithfully} on $\cD$.

If the group $G$ acts on the design $\cD$ as a group of automorphisms
then we get three induced permutation groups, namely:

\begin{enumerate}
\item $(G, P_\cD)$ -- The induced permutation group on the points.
\item $(G, B_\cD)$ -- The induced permutation group on the blocks.
\item $(G, F_\cD)$ -- The induced permutation group on the flags.
\end{enumerate}

We will be interested in highly symmetric designs.
In particular, we will be interested in pairs $(G, \cD)$
where $G$ acts on $\cD$ in such a way that the induced
permutation group on the flags $(G, F_\cD)$ is transitive.
We will either call the pair $(G, \cD)$ a flag-transitive design,
or we call the design $\cD$ a $G$-flag transitive design.
Many interesting groups arise naturally as flag-transitive
automorphism groups of designs, including a number
of the sporadic simple groups.

\begin{proposition}
If $(G,\cD)$ is a flag-transitive design
then the induced permutation groups $(G,P_\cD)$ and $(G,B_\cD)$
are both transitive
\end{proposition}

%Suppose that $\cD$ is a $G$-flag transitive design, and $(p,b)$ is some flag.
%If $H$ is the stabilizer of $p$, then $(G,P_\cD)$ is a transitive permutation group
%which is permutationally isomorphic to $(G,\cos(H))$.
%If $S$ is the stabilizer of $b$, then $(G,B_\cD)$ is a transitive permutation group
%which is permutationally isomorphic to $(G,\cos(S))$.

\begin{proposition}
If $(G,\cD)$ is a flag-transitive design and $H_1 < G$
is the stabilizer of some point $p$ then the
induced permutation group $(H_1, T(p))$ is transitive.
Simmilarly if $H_2 < G$ is the stabilizer of some block $b$
then the induced permutation group $(H_2, T(b))$ is transitive.
\end{proposition}

%Let $R = S \cap H$, then $R$ stabilizes the flag $(p,b)$.
%Let $B(p)$ denote the blocks of $\cD$ incident with $p$
%and let $P(b)$ denote the points of $\cD$ incident with $b$.
%The action of $H$ on $B(p)$ is permutationally isomorphic to $(H,\cos(R))$
%and the action of $S$ on $P(b)$ is permutationally isomorphic to $(S, \cos(R))$.

If $\cD_1$ and $\cD_2$ are two $G$-flag transitive designs, 
then a {\it $G$-design homomorphism} from $\cD_1$ to $\cD_2$
is a pair of maps $\rho = (\rho_P, \rho_B)$
such that:

\begin {enumerate}
\item $\rho_P : P_{\cD_1} \to P_{\cD_2}$ induces a permutation homomorphism between 
$(G,P_{\cD_1})$ and $(G,P_{\cD_2})$.
\item $\rho_B : B_{\cD_1} \to B_{\cD_2}$ induces a permutation homomorphism between 
$(G,B_{\cD_1})$ and $(G,B_{\cD_2})$.
\item $(p,b) \in I_1$ if and only if $(\rho_P(p), \rho_B(b)) \in I_2$.
\end{enumerate}

When both $\rho_P$ and $\rho_B$ are bijections we have a {\it $G$-design isomorphisms}.
When confusion is unlikely to result, 
we will sometimes write $\rho(p)$ instead of $\rho_P(p)$ when $p$ is a point
and $\rho(b)$ instead of $\rho_B(b)$ when $b$ is a block.

\begin{proposition}
If $\rho$ is a $G$-design homomorphism from $\cD_1$ to $\cD_2$
then $\rho$ induces a $G$-permutation homomorphism from $F_{\cD_1}$ to $F_{\cD_2}$.
\end{proposition}

%For a fixed $G$ the set of all $G$-flag transitive designs
%with $G$-design homomorphisms between them form a category which we shall denote by
%$G$-Flag.
\section{$G$-symmetric designs}

In this section we show how $G$-symmetric designs
may be viewed as a special kind of flag-transitive design.
I'm calling them symmetric designs, but I think
they should probably be called square designs.

\begin{definition}
A $G$-flag transitive design $\cD$ is {\it self-dual}
if there exists a $G$-isomorphism $\rho = (\rho_P, \rho_B)$
between $\cD$ and $\cD^{*}$.
The $G$-isomorphism $\rho$ is called a {\it duality} of $\cD$.
\end{definition}

\begin{definition}
A {\it $G$-symmetric design} is self-dual $G$-flag transitive design $\cD$
which admits a dulaity $\rho = (\rho_P, \rho_B)$ with the property that
$\rho_B \circ \rho_P = id_P$ and $\rho_P \circ \rho_B = id_B$.
In this case the $G$-isomorphism $\rho$ is called a {\it polarity}
of $\cD$.
\end{definition}

Observe that if $\cD$ is a $G$-symmetric design with the property
that $(p, \rho(p))$ is a a flag for some $p \in P_\cD$,
then since $\rho$ is a $G$-isomorphism it follows
that the point-stabilizer of $\cD$ is isomorphic to the flag-stabilizer of $\cD$.
Thus each point is incident with exactly one block
and vice versa.
We consider such $G$-symmetric designs to be ``degenerate'',
and unless an explicit statement to the contrary is given shall
take the expression ``$G$-symmetric design''
to mean ``non-degenerate $G$-symmetric design''.

In a sense that will become clear a little later,
these ``degenerate'' $G$-symmetric
designs correspond to the ``degenerate'' orbital graph $\orb_\Delta(G,\Omega)$
which could be formed from the permutation group $(G, \Omega)$
by taking $\Delta$ to be the diagonal orbit,
and also to the ``degenerate'' Sabidussi graph $\sab(G,H,a)$
which could be formed by taking $a \in H$.
\hfill \break

We shall see that each $G$-symmetric graph
gives rise in a natural way to a $G$-symmetric design
and conversely a $G$-symmetric design together
with a ``marked'' polarity give rise to a $G$-symmetric graph.
We shall also see that in some cases, by choosing a different
polarity it is possible to construct two non-isomorphic $G$-symmetric
graphs from the same $G$-symmetric design.
For any graph $\Gamma$, let $N \Gamma = \{ \Gamma(v) : v \in V \Gamma \}$
denote the set of {\it neighbourhoods} of $\Gamma$.

\begin{proposition}
If \/ $\Gamma$ is a $G$-symmetric graph,
then the incidence structure $\cD(\Gamma) = (V \Gamma, N \Gamma,I)$
where $(v,n) \in I$ if and only if $v \in n$ is a $G$-symmetric design
\begin{proof}
The $G$-arc transitivity of $\Gamma$ is sufficient to ensure
that $\cD(\Gamma)$ is a $G$-flag transitive design.
To see that it is in fact
a $G$-symmetric design we must exhibit a polarity.

Let $\rho_P : P \to B$ be given by $v \mapsto \Gamma(v)$.
The map $\rho_P$ is clearly bijective, and since $G$ acts on $\Gamma$
as a group of automorphisms we have:
\begin{eqnarray*}
\rho_P(v^g) & = & \Gamma(v^g) \\
            & = & \Gamma(v)^g \\
	                & = & \rho_P(v)^g.
			\end{eqnarray*}
			So $\rho_P$ induces a permutation isomorphism.
			Now, take $\rho_B = \rho_P^{-1}$.
			We must check that the pair $\rho = (\rho_P, \rho_B)$ preserve
			the incidence structure of the design.

Suppose that $(v,\Gamma(w))$ is a flag of $\cD$, so $v \in \Gamma(w)$.
Since $\Gamma$ is a simple graph it follows immediately that $w \in \Gamma(v)$
and so $(w, \Gamma(v))$ is also a flag of $\cD$.
That is $(\Gamma(v), w)$ is a flag of $\cD^{*}$.
But
\begin{eqnarray*}
\rho((v, \Gamma(w))) & = & (\rho_P(v), \rho_B(\Gamma(w))) \\
                     & = & (\Gamma(v), w).
		     \end{eqnarray*}
		     So we are done.

		     \end{proof}
		     \end{proposition}

\begin{proposition}
If $\cD$ is a (non-degenerate) $G$-symmetric design with ``marked'' polarity $\rho$
then the graph $\Gamma(\cD, \rho)$ with vertex set $P_\cD$
and arc set $\{ (p,q) : (q,\rho_P(p)) \in I_\cD \}$ is $G$-symmetric.
\begin{proof}
We must first check that $\Gamma(\cD, \rho)$ is well-defined.
The non-degenerateness of $\cD$ ensures that there are no loops.
If $(p,q)$ is an arc of $\Gamma(\cD, \rho)$
then $(q, \rho_P(p))$ is a flag of $\cD$.
Since $\rho$ is an isomorphism,
if $(q, \rho_P(p))$ is a flag of $\cD$
then $(\rho_P(q), \rho_B \circ \rho_P(p)) = (\rho(q), p)$
is a flag of $\cD^{*}$.
It follows that $(p,\rho(q))$ is a flag of $\cD$,
and thus $(q,p)$ is an arc of $\Gamma(\cD, \rho)$.
So $\Gamma(\cD, \rho)$ is simple.

By proposition [?] $G$ acts transitively on the points $\cD$
so $\Gamma(\cD, \rho)$ is $G$-vertex transitive.
Suppose that $(p,q_1)$ and $(p,q_2)$ are two distinct arcs of
$\Gamma(\cD, \rho)$. Then we have $q_1, q_2 \in T(\rho(p))$.
By proposition [?]
the stabilizer of $\rho(p)$ acts transitively on $T(\rho(p))$.
Thus we can find some $g \in G$ which fixes $\rho(p)$ and
carries $q_1$ to $q_2$.
Since $\rho$ is a $G$-isomorphism, if $g$ fixes $\rho(p)$
it also fixes $p$, thus $g$ carries the arc $(p,q_1)$
to the arc $(p,q_2)$. So $\Gamma(\cD, \rho)$ is $G$-locally
transitive. The result follows.

\end{proof}
\end{proposition}

\begin{proposition}
For any design $\cD$ and any polarity $\rho$ the design $\cD(\Gamma(\cD,\rho))$
is isomorphic to $\cD$
\begin{proof}
later
\end{proof}
\end{proposition}

%\begin{definition}
%If $\cD_1$ and $\cD_2$ are two $G$-symmetric designs
%then a $G$-design homomorphism from $\cD_1$ to $\cD_2$
%is a pair of permutation homomorphisms $\eta_P : P_1 \to P_2$ and $\eta_B : B_1 \to B_2$
%satisfying $(p,b) \in I_1$ implies $(\eta_P(p),\eta_B(b)) \in I_2$.
%\end{definition}

%\begin{proposition}
%If $\cD_1$ and $\cD_2$ are two $G$-symmetric designs
%then a $G$-design homomorphism from $\cD_1$ to $\cD_2$
%induces a a $G$-permutation homomorphism from the flags of $\cD_1$
%to the flags of $\cD_2$
%\end{proposition}

%What I want to say is that $G$-symmetric designs
%together with $G$-design homomorphisms form a category
%(well, kind of more of a poset really).
%I also want to say that the category $G$-Graph
%``projects down'' onto the category $G$-Design.
%$F : G$-Graph $\to G$-Design.
%If there's a $G$-graph homomorphism from $x$ to $y$
%then there's a $G$-design homomorphism from $F(x)$ to $F(y)$

\section{``cross-section'' of a graph homomorphism}

Suppose that $\Gamma$ is an imprimitive $G$-symmetric graph
with $\cB$ a non-trivial $G$-invariant partition of the vertices.
Assume further that the quotient $\Gamma_\cB$ has valency
at least one so that the blocks of $\cB$ are independent sets.

For any vertex $\alpha$ of $\Gamma$, let $B(\alpha)$ denote
the block of $\cB$ containing $\alpha$.
Let $\Gamma(\alpha) = \{ \beta \in V \Gamma : (\alpha, \beta) \in A \Gamma \}$
denote the neighbours of $\alpha$ in $\Gamma$.
Let $\Gamma_\cB(B) = \{ C \in \cB : (B,C) \in A \Gamma_\cB \}$
denote the neighbours of $B$ in the quotient.
For any $\alpha \in B$, let
$\Gamma_\cB(\alpha) = \{ B(\beta) : \beta \in \Gamma(\alpha) \}$
denote the blocks of $\cB$ containing the neighbours of $\alpha$ in $\Gamma$.

Construct a design with point set $B$, block set $\Gamma_\cB(B)$
and an incidence relation $ I \subseteq B \times \Gamma_\cB(B)$
defined by $(\alpha, C) \in I$ if and only if $C \in \Gamma_\cB(\alpha)$.

\begin{proposition}
The incidence structure $\cD(B) = (B, \Gamma_\cB(B), I)$ is a $H$-flag transitive
$(v,k,\lambda)$ design, where:
\begin{eqnarray*}
v & := &  |B| \\
k & := &  |\Gamma(C) \cap B|, \\
\lambda  & := &  |\Gamma_\cB(\alpha)| \\
\end{eqnarray*}
and $H$ is the stabilizer of the block $B \in \cB$.
\end{proposition}
	      
The design $\cD(B)$ gives, in a sense, a ``cross-section''
of the homomorphism $\Gamma \mapsto \widetilde{\Gamma}$.

\section{Reconstruction problem}

Gardiner and Praeger observed that if $\Gamma$ is a $G$-symetric graph
and $\cB$ is a non-trivial $G$-invariant partition of the vertices,
then the graph $\Gamma$ ``decomposes''
into the triple $(\Gamma_\cB, \Gamma[B,C], \cD(B))$.

We may ask now, to what extent the combinatorial structure of $\Gamma$
is determined by the triple $(\Gamma_\cB, \Gamma[B,C], \cD(B))$?
Do these three components contain sufficient information
to reconstruct $\Gamma$ uniquely?
If not, what extra information is required?

%The problem then becomes that of reconstructing $\Gamma$ 
%from the triple $(\Gamma_\cB, \Gamma[B,C], \cD(B))$.

Suppose we are given a symmetric graph $\Lambda$,
a symmetric bipartite graph $\Sigma$ and a flag transitive design
$\cD$ without any particular groups acting upon them.
Let us say that a symmetric graph $\Gamma$ is 
``product'' of $(\Lambda, \Sigma, \cD)$ if there is a group $G$
acting $\Gamma$ and a $G$-invariant partition $\cB$ of
the vertices of $\Gamma$
such that the triple $\Gamma$ decomposes into the triple
$(\Lambda, \Sigma, \cD)$ 

What are the necessary and sufficient
conditions under which these three components $(\Lambda, \Sigma, \cD)$
can be ``glued together'' into some ``product'' $\Gamma$ ?

If the necessary conditions are satisfied,
could there be more than one way to ``glue together'' a given 
triple $(\Lambda, \Sigma, \cD)$ into an imprimitive symmetric graph 
$\Gamma$ ?

%% file: threearc.tex
\chapter{Some Constructions}

% Try to weave references to the literature into the narrative
The problem of relating the structure of $\Gamma$ 
to that of the triple $(\Gamma_\cB, \Gamma[B,C], \cD(B))$
is very difficult when $\Gamma$ is taken to be 
an arbitrary imprimitive symmetric graph.
The approach which has been taken by researchers in the area
is to first impose additional assumptions 
on one or more of the components of the decomposition,
then study the subfamily of imprimitive symmetric graphs 
admitting a decompositions which satisfies these additional assumptions.

The special case where the parameters of the ``kernel design'' $\cD(B)$
satisfy: $$v = k-1 \geq 2.$$
was studied by Li, Praeger and Zhou \cite{ZCP2000}.
It was found that for this special case, 
if the design $\cD(B)$ contains no repeated blocks,
then there exists an elegant combinatorial method
for constructing $\Gamma$ from $\Gamma_\cB$.
This is the {\it three-arc graph construction} which will
be described in the next section.

In a later paper Zhou \cite{ZDM2002} showed that 
the three-arc graph construction actually applies 
to a wider family of triples $(G, \Gamma, \cB)$
than those originally studied.
In particular the construction may be used for any triple $(G, \Gamma, \cB)$
satisfying the following condition:  

\begin{condition}
The induced actions of $G_B$ on $B$ and $\Gamma_\cB(B)$
are permutationally equivalent with respect to some bijection
$\rho : B \to \Gamma_\cB(B)$.
\end{condition}

On a group theoretic level this is quite a natural condition to impose.
Suppose that the triple $(G,\Gamma, \cB)$ satisfies the above condition
and that $\Gamma \cong \Sab(G,K,a)$ 
and $\Gamma_\cB \cong \Sab(G,H,a)$.
If $B$ is a block of $\cB$ then
by proposition [?] the action of $G_B$ on $B$
is permutation equivalent to $(H, \cos(K))$.
The action of $G_B$ on $\Gamma_\cB(B)$ 
is just the ``local action'' of $\Gamma_\cB$
and so the permutation group $(G_B, \Gamma_\cB(B))$
is permutation equivalent to $(H, \cos(a^{-1}Ha \cap H))$.

If the condition above holds, 
then we must have $K \cong a^{-1}Ha \cap H$.
This means that the stabilizer of a point of $\Gamma$ is isomorphic
to the stabilizer of an arc of $\Gamma_\cB$.
This strong ``coupling'' between the points of $\Gamma$
and the arcs $\Gamma_\cB$ allows us to construct
$\Gamma$ from $\Gamma_\cB$ in a straightforeward manner.

%The condition above is also very natural on a design theoretic level.
%Recall that if $\Gamma \cong \Sab(G,K,a)$ and $\Gamma_\cB \cong \Sab(G,H,a)$
%then the ``kernel design'' of the homomorphism $\Gamma \mapsto \Gamma_\cB$
%has point stabilizer $a^{-1}Ha \cap H$ and block stabilizer $\cK$.
%If $\cK \cong a^{-1}Ha \cap H$ then the action of $H$ on the points
%of $\cK$ is permutation equivalent to the action of $H$ on the blocks
%of $\cK$.

\section{Three-Arc Graphs}
The three-arc graph construction was first introduced 
by Li, Praeger and Zhou in \cite{ZCP2000}.
Let $\Sigma$ be any $G$-symmetric graph.
Recall that an {\it $s$-arc} of $\Sigma$ is a sequence 
$(\alpha_0, \alpha_1, \dots, \alpha_s)$ of vertices in $\Sigma$ such that 
$\alpha_i,\alpha_{i+1}$ are adjacent in $\Sigma$ and $\alpha_{i-1} \neq \alpha_{i+1}$
for each $i$. The set of $s$-arcs of $\Sigma$ is denoted by $\Arc_s(\Sigma)$.
Consider the induced action of $G$ on $\Arc_3(\Sigma)$ given by:
$$(\sigma_0,\sigma_1,\sigma_2,\sigma_3)^g = (\sigma_0^g,\sigma_1^g,\sigma_2^g,\sigma_4^g).$$
This action is, in general, intransitive. We may however partition the set $\Arc_3(\Sigma)$
into {\it orbits} on which $G$ acts transitively. For such an orbit $\Delta$, let 
$\Delta^{\circ} = \{ (\sigma_4, \sigma_3, \sigma_2, \sigma_1) : 
                 (\sigma_1, \sigma_2, \sigma_3, \sigma_4) \in \Delta \}$ 
denote the {\it pair} of $\Delta$. It is not hard to check that
$\Delta^{\circ}$ is again an orbit of $G$ on $\Arc_3(\Sigma)$.
If $\Delta = \Delta^{\circ}$ then $\Delta$ is said to be {\it self paired}.

\begin{definition}
\label{threearc}
Given a $G$-symmetric graph $\Sigma$ and a self-paired orbit $\Delta$ on 
$\Arc_3(\Sigma)$, the {\it three-arc graph}
$\Gamma = \Arc_{\Delta}(\Sigma)$ is the graph with vertex set $A \Sigma$
and arc set $\{ ((\sigma,\tau), (\sigma',\tau') : (\tau,\sigma,\sigma',\tau') \in \Delta \}$.
\end{definition}

Note that the requirement that $\Delta$ is self-paired 
ensures that the resulting graph is simple.

\begin{proposition}
With $G$, $\Sigma$ and $\Delta$ as in definition \ref{threearc},
the three-arc graph $\Gamma = \Arc_{\Delta}(\Sigma)$
is $G$-symmetric.
\begin{proof} Immediate from the construction.
\end{proof} 
\end{proposition}

For each vertex $\sigma$ of $\Sigma$, let 
$B(\sigma) = \{ (\sigma, \tau) : \tau \mbox{ is a neighbour of } \sigma \}$
be the set of arcs of $\Sigma$ with initial vertex $\sigma$.
Clearly $\cB = \{ B(\sigma) : \sigma \mbox{ is a vertex of } \Sigma \}$ 
is a partition of $A \Sigma$.  
For any $g \in G$ we have 
\begin{eqnarray*}
B(\sigma)^g & = & \{ (\sigma, \tau)^g : \tau \mbox{ is a neighbour of } \sigma \} \\
              & = & \{ (\sigma^g, \tau^g) : \tau \mbox{ is a neighbour of } \sigma \} \\
	      & = & \{ (\sigma^g, \tau) : \tau \mbox{ is a neighbour of } \sigma^g \} \\
	      & = & B(\sigma^g).
\end{eqnarray*}
Thus $\cB = \{ B(\sigma) : \sigma \in V \Sigma \}$ 
is in fact a $G$-invariant partition of $V \Gamma$. 

\begin{proposition}
\label{three-arc quotient}
With $\Gamma = \Arc_{\Delta}(\Sigma)$ and $\cB$ as defined above,
The quotient graph $\Gamma_\cB$ is isomorphic to $\Sigma$.

\begin{proof}
The map $\sigma \mapsto B(\sigma)$ identifies the vertices 
of $\Sigma$ with the vertices of $\Gamma_\cB$.
We must show that 
$(\sigma, \sigma')$ is an arc of $\Sigma$ if and only if
$(B(\sigma), B(\sigma'))$ is an arc of $\Gamma_\cB$.

If $(B(\sigma), B(\sigma'))$ is an arc of $\Gamma_\cB$
then then for some $\tau, \tau' \in V \Sigma$ 
the arcs $(\sigma, \tau)$ and $(\sigma',\tau')$ are adjacent in $\Gamma$.
It follows that 
$(\tau, \sigma, \sigma',\tau')$ is a three-arc of $\Sigma$
and so $\sigma$ is adjacent to $\sigma'$ in $\Sigma$. 

For the other direction, 
let $(\alpha, \beta, \gamma, \delta)$ be any three-arc 
of $\Sigma$ contained in in $\Delta$.
If $(\sigma, \sigma')$ is any arc of $\Sigma$, 
then by the $G$-arc-transitivity of $\Sigma$
there exists a $g \in G$ such that $(\beta,\gamma)^g = (\sigma,\sigma')$.
Since $\Delta$ is an orbit of $G$ on the three-arcs of $\Sigma$, 
it follows that 
$(\alpha,\beta,\gamma,\delta)^g = (\alpha^g, \sigma, \sigma', \delta^g)$ 
is contained in $\Delta$.
Thus the arc $(\sigma, \alpha^g)$ is adjacent to 
the arc $(\sigma', \delta^g)$ in $\Gamma$
and so $B(\sigma)$ is adjacent to $B(\sigma')$ in $\Gamma_\cB$.

\end{proof}
\end{proposition}

%The three-arc graph construction allows us to take a $G$-symmetric graph
%$\Sigma$ and construct a {\it larger} imprimitive $G$-symmetric graph $\Gamma$
%which admits the original graph $\Sigma$ as a quotient.
%Provided of course that we can find a suitable self-paired orbit
%on the three arcs of $\Sigma$.
%This is an interesting result because blah blah blah....

%We shall now show that the family of triples $(G, \Gamma, \cB)$
%for which $\Gamma$ is a three-arc graph of $\Gamma_\cB$ 
%is precisely the family of triples satisfying 
%the condition {\bf (PE)} stated earlier.

\begin{proposition}

If $(G,\Gamma, \cB)$ is such that $\Gamma$ is a three-arc graph
of $\Gamma_\cB$, then for any $B \in \cB$
the permutation groups $(G, B)$ and $(G, \Gamma_\cB(B))$
are permutation equivalent.  

\begin{proof}
Each $B \in \cB$ is of the form $B(\sigma)$ for some
$\sigma$ a vertex of $\Gamma_\cB$.
The vertices of $\Gamma$ contained in $B(\sigma)$ 
are precisely the arcs of $\Gamma_\cB$
of the form $(\sigma, \tau)$.

Let $\rho : B(\sigma) \to \Gamma_\cB(B(\sigma))$ be given by 
$(\sigma,\tau) \mapsto B(\tau)$.
Clearly $\eta$ is a bijection.
If $g$ is any element of the stabilizer of $B(\sigma)$ 
then we have:
\begin{eqnarray*}
\rho((\sigma,\tau)^g) & = & \rho((\sigma^g,\tau^g)) \\
                      & = & \rho((\sigma,\tau^g))   \\
		      & = & B(\tau^g)               \\
		      & = & B(\tau)^g.
\end{eqnarray*}
Thus $\rho$ induces a permutation equivalence between the permutation groups
$(G_{B(\sigma)}, B(\sigma))$ and $(G_{B(\sigma)}, \Gamma_\cB(B(\sigma))$.

\end{proof}
\end{proposition}

\begin{proposition}
If $(G, \Gamma, \cB)$ is such that 
$(G, B)$ is permutation equivalent to $(G, \Gamma_\cB(B))$
for some $B \in \cB$, then the action of $G$
on the vertices of $\Gamma$ is permutation equivalent
to the action of $G$ on the arcs of $\Gamma_\cB$.

\begin{proof}
Let $\rho : B \to \Gamma_\cB(B)$ be any bijection
inducing a permutation equivalence between 
$(G,B)$ and $(G, \Gamma_\cB(B))$
and let $\pi : V \Gamma \to \cB$ be the map which sends
each vertex of $\Gamma$ to the block of $\cB$ containing it.
Fix a vertex $\alpha$ of $\Gamma$.
Since $\Gamma$ is $G$-vertex transitive, every vertex of $\Gamma$
may be written in the form $\alpha^g$ for some $g \in G$.
Let $\mu : V \Gamma \to A \Gamma_{\cB}$ be given by
$ \alpha^g \mapsto (\pi(\alpha)^g, \rho(\alpha)^g) $.
It is not too hard to see that $\mu$ is a bijection.
For any $\beta = \alpha^g$ and any $x \in G$ 
we have:
\begin{eqnarray*}
\mu(\beta^x) & = & \mu(\alpha^{gx})                       \\
             & = & (\pi(\alpha)^{gx}, \eta(\alpha)^{gx})    \\
	     & = & (\pi(\alpha^g)^x, \eta(\alpha^g)^x)      \\
	     & = & (\pi(\beta)^x, \eta(\beta)^x)            \\
	     & = & (\pi(\beta), \eta(\beta))^x              \\
	     & = & \mu(\beta)^x
\end{eqnarray*}
Thus $\mu$ establizhes a permutation isomorphism between
$(G,V \Gamma)$ and $(G, A \Gamma_\cB)$.

\end{proof}
\end{proposition}

We may make use the map $\mu$ to ``label'' the vertices
of $\Gamma$ by the arcs of $\Gamma_\cB$.
For any arc $(B,C)$ of $\Gamma_{\cB}$, 
let $v_{BC} = \mu^{-1}((B,C))$ denote
the vertex of $\Gamma$ which is mapped to $(B,C)$ by $\mu$.

\begin{proposition}
\label{labelling}
Some stuff about how $G$ acts on labelled vertices
and how the initial block contains the vertex.
\end{proposition}

\begin{proposition}
\label{PE implies 3-arc}
Provided that $\Gamma$ as valency at least two,
for each arc $(v_{BC}, v_{DE})$ of $\Gamma$,
$(C,B,D,E)$ is a 3-arc of $\Gamma_\cB$.

\begin{proof}
Suppose that $v_{BC}$ were adjacent to $v_{CB}$.
Since $\val(\Gamma) \geq 2$, the vertex $v_{CB}$
is adjacent to some other vertex $v_{B_1C_1}$ 
distinct from $v_{BC}$
By the $G$-symmetry of $\Gamma$ there exists a $g \in G$
such that $(v_{BC},v_{CB})^g = (v_{B_1C_1},v_{CB})$.
By the previous proposition, this implies that 
$B = B^x = B_1$ and $C = C^x = C_1$, so $v_{BC} = v_{B_1C_1}$.
A contradiction. Thus $v_{BC}$ is not adjacent to $v_{CB}$.

Suppose now that $v_{BC}$ were adjacent to $v_{CE}$
for some $E \neq B$.
Since $\val(\Gamma) \geq 2$ we can find another vertex $v_{C_1E_1}$
distinct from $v_{CE}$ and adjacent to $v_{BC}$.
By the $G$-symmetry of $\Gamma$ we can find some $g_1 \in G$
such that $(v_{BC},v_{CE})^{g_1} = (v_{BC},v_{C_1E_1})$.
By the previous proposition, this implies that $C = C^{g_1} = C_1$.
We can also find a $g_2 \in G$
such that $(v_{BC},v_{CE})^{g_2} = (v_{C_1E_1}, v_{BC})$
and so 
$B = C^{g_2} = E_1$. Combining these gives $v_{C_1E_1} = v_{CB}$, 
but $v_{BC}$ is not adjacent to $v_{CB}$
so again we have a contradiction.

We know now that if $(v_{BC}, v_{DE})$ is an arc of $\Gamma$,
then $D \neq C$. A simmilar argument shows that $B \neq E$,
so $B,C,D,E$ are distinct vertices with $B$ adjacent to $C$
and $D$ adjacent to $E$. Since, by the previous proposition
$B(v_{BC}) = B$ and $B(v_{DE} = D$, we have $B$ adjacent to $D$
and so $C,B,D,E$ is a 3-arc of $\Gamma_\cB$ as claimed.

\end{proof}
\end{proposition}

Concluding remarks.

%It follows immediately from proposition \ref{PE implies 3-arc}
%that if $(\Gamma, \cF)$ is any pair satisfying {\bf PE assumption}
%then $\Gamma$ is isomorphic to a 3-arc graph of $\Gamma_\cF$
%taking $\Delta$ to be the orbit containing $(C,B,D,E)$
%for any arc $(v_{BC}, v_{DE})$ of $\Gamma$.
%(This needs to be much better written).
%Need to talk about conditions under which a $G$-symmetric graph
%admits a self-paired $G$-orbit on its three-arcs.

%\section{$(G,2)$-arc transitive quotients}

%\section{$(G,3)$-arc transitive quotients}

%\section{Almost Covers}

%\section{The case when the ``kernel design'' is symmetric}

%\section{Uniqueness of the Construction}

%\begin{enumerate}
%\item how to find the self-paired orbit?
%\item not unique
%\item extremal case: induced bipartite graph almost a matching
%\item occurs only if quotient is $(G,3)$-arc transitive
%\item which paper does this occur in?
%\end{enumerate}

%We have seen that any pair $(\Gamma, \cB)$ satisfying the
%{\it (PE) assumption}, $\Gamma$ may be reconstructed
%from $\Gamma_\cB$ by using the three-arc graph construction
%{\it for some suitable self-paired orbit}. 
%The question is, which orbit? and is there a unique such orbit?

%An obvious case in which there will be a unique orbit
%on the three arcs is when the quotient is in fact
%$(G,3)$-arc transitive. In this case, $G$ is in fact
%transitive on $\Arc_3(\Gamma_\cB)$ and the set
%of all $3$-arcs forms a self-paired orbit. 
\section{Locally imprimitive quotient}

Suppose that $\Gamma = \sab(G,H,a)$.
Let $\overline{H} = a^{-1}Ha \cap H$,
so the local permutation group induced at each vertex of $\Gamma$
is equivalent to $(H, \cos(\overline{H}))$.
Suppose further that there is some $K$ such that $\overline{H} < K < H$
and $a \not\in K$. Then the graph $\widetilde{\Gamma} = \sab(G,K,a)$
is an extension of $\Gamma$.
Let $\overline{K} = a^{-1}Ka \cap K$.
Since $K < H$ we must also have $a^{-1}Ka < a^{-1}Ha$
so $\overline{K} \leq \overline{H}$.
Since $\overline{H} < K$ we must also have
$a^{-1}\overline{H}a < a^{-1}Ka$
but since $a$ is an involution $a^{-1}\overline{H}a = \overline{H}$
so we have $\overline{H} \leq \overline{K}$.
That is $\overline{H} = \overline{K}$
and the local permutation group induced at each vertex of
$\widetilde{\Gamma}$ is $(K, \cos(\overline{H}))$.

The pair $(\Gamma,\widetilde{\Gamma})$ satisfy the property that
{\it globally} $\widetilde{\Gamma}$ is imprimitive,
admitting $\Gamma$ as quotient, but {\it locally}
$\Gamma$ is imprimitive, and the local permutation group of
$\widetilde{\Gamma}$ is a quotient of the local permutation
group of $\Gamma$. We wish to understand how $\Gamma$
and $\widetilde{\Gamma}$ are related {\it structurally},
and ideally find some method of constructing $\widetilde{\Gamma}$
from $\Gamma$.

As a preliminary observation, let $n = [G:H]$ be the number of vertices of $\Gamma$,
let $v = [K:\overline{H}]$ be the valency of $\widetilde{\Gamma}$
and let $r = [H:K]$.
The graph $\widetilde{\Gamma}$ has $[G:K] = [G:H][H:K] = nr$
vertices. That is, $r$ times as many vertices as $\Gamma$.
The valency of $\Gamma$ is $[H:\overline{H}] = [H:K][K:\overline{H}] = rv$.
That is $r$ times the valency of $\widetilde{\Gamma}$.
Since the number of edges in a graph is equal to half the
valency times the number of vertices, it follows that
both $\Gamma$ and $\widetilde{\Gamma}$ have the same number of edges.
\hfill \break

The action of $G$ on the arcs of $\Gamma$ is permutation equivalent
to $(G, \cos(\overline{H}))$.
Since $\Gamma$ is simple, each arc $(\alpha,\beta)$ has a {\it pair}
namely $(\beta,\alpha)$.
The function $\varphi : A \Gamma \to A \Gamma$
which sends each arc to its pair is an involution
It also preserves the action of $G$.

Now, since $\overline{H} < H < G$ the action of $G$
on the arcs of $\Gamma$ is imprimitive.
That is, there is some $G$-invariant partition $\cB$
of $A \Gamma$ such that the action of $G$ on $\cB$
is permutation equivalent to $(G, \cos(H))$.
In fact, this is the partition $\cB = \{ B(\alpha) : \alpha \in V \Gamma \}$
where $B(\alpha)$ is the set of all arcs who's initial vertex is $\alpha$.

Let $\pi : A \Gamma \to \cB$ be the map which sends
each arc to the block of $\cB$ containing it.
We may define a new graph $\Sigma$ we vertex set
$\cB$ where $p$ is adjacent to $b$ if and only if
there is some $x \in \pi^{-1}(p)$ and some $y \in \pi^{-1}(b)$
such that $y = \varphi(x)$.
This new graph $\Sigma$ is in fact isomorphic to $\Gamma$.

The action of $G$ on the arcs of $\Gamma$ is permutation
isomorphic to the action of $G$ on the arcs of $\widetilde{\Gamma}$.
Both actions are of the form $(G, \cos(\overline{H}))$
Since $\overline{H} < K < G$ there exists some $G$-invariant partition
$\mathfrak{B}$ of $A \Gamma$ such that the action of $G$ on $\mathfrak{B}$
is permutation equivalent to $(G, \cos(H))$.
The partition $\mathfrak{B}$ is a refinement of the partition $\cB$.

Let $\lambda : A \Gamma \to \mathfrak{B}$ be the map hich sends
each arc to the block of $\mathfrak{B}$ containing it.
We may define a new graph $\widetilde{\Sigma}$ we vertex set
$\mathfrak{B}$ where $p$ is adjacent to $b$ if and only if
there is some $x \in \lambda^{-1}(p)$ and some $y \in \lambda^{-1}(b)$
such that $y = \varphi(x)$.
This new graph $\Sigma$ is in fact isomorphic to $\widetilde{\Gamma}$.
Thus we have a method of construction $\widetilde{\Gamma}$ from $\Gamma$.

\section{Covering Graphs}

% Introductory paragraph:
The three-arc graph construction allows us, under certain conditions,
to ``unfold'' a $G$-symmetric graph into a larger, {\it imprimitive}
$G$-symmetric graph, admitting the original graph as a quotient.
The {\it covering graph} construction \cite{B1974} is simmilar
in this respect. It also, under certain conditions,
allows us to ``unfold'' a $G$-symmetric graph into a larger one.
It differs from the three-arc graph construction, however,
in that it requires a simultaneous ``unfolding''
of the group $G$.

% Semidirect Product:
Recall that given two groups $N$ and $G$ and a homomorphism
$$\rho : G \to \Aut(N)$$
we may form the {\it semidirect product} of $N$ by $G$.
This is the group $\widetilde{G} = N \rtimes_\rho G$,
whose elements are the ordered pairs:
$$\{ (n,g) : n \in N, g \in G \}$$
with multiplication given by:
$$(n_1, g_1)(n_2,g_2) = (n_1^{\rho(g_1)}n_2, g_1g_2).$$
The functions $i_1 : N \to \widetilde{G}$ and $i_2 : G \to \widetilde{G}$
given by $ n \mapsto (n,1) $ and $ g \mapsto (1,g) $ respectively give
natural embeddings of $N$ and $G$ into $\widetilde{G}$.
Identifying $N$ with its image under $i_1$,
the semidirect product has the property that $N$ is normal in $\widetilde{G}$
and $\widetilde{G}/N \cong G$.

% Definitions:
\begin{definition} [$N$-chain]
Suppose that $\Gamma$ is a $G$-symmetric graph and $N$ is a group.
An $N$-chain is a function $\phi : A \Gamma \to N$
satisfying $\phi((v,u)) = \phi((u,v))^{-1}$ for all $(u,v) \in A \Gamma$.
\end{definition}

\begin{definition} [Compatible $N$-chain]
Suppose that $\Gamma$ is a $G$-symmetric graph, $N$ is a group
and $\rho : G \to \Aut(N)$ is a homomorphism.
An $N$-chain $\phi$ is said to be {\it compatible} with $\rho$
if for every $g \in G$ the following diagram commutes:

\centerline{
\xymatrix{
\ar[d]_{g} \ar[r]^\phi A \Sigma & K \ar[d]^{\rho(g)} \\
\ar[r]_{\phi} A \Sigma & K
}}
\end{definition}

\begin{definition}[Biggs Cover]
Suppose that $\Gamma$ is a $G$-symmetric graph, $N$ is a group,
$\rho : G \to \Aut(N)$ is a homomorphism and $\phi$ is a compatible $N$-chain.
The Biggs Cover $\widetilde{\Gamma}(N,\rho,\phi)$
of $\Gamma$ with respect to $\phi$ is the graph
with vertex set:
$\{ (g,v) : k \in G, v \in V \Gamma \}$
and arc set:
$ \{ ((g_1,v_1),(g_2,v_2)) : (v_1,v_2) \in A \Gamma, g_2 = g_1 \phi((v_1,v_2)) \}$
\end{definition}

Note that the condition $\phi((v,u)) = \phi((u,v))^{-1}$ ensures that the resulting
graph is simple.

% Proposition:
\begin{proposition}
The covering graph $\widetilde{\Gamma}(N,\rho,\phi)$ is $(N \rtimes_\rho G)$-symmetric.
\begin{proof}

Define an action of $N \rtimes_\rho G$ on $\widetilde{\Gamma}$ by:
$$(n,v)^{(\eta, g)} = (n^{\rho(g)} \eta,v^g).$$
This action is well-defined since:
\begin{eqnarray*}
(n,v)^{(\eta_1, g_1)(\eta_2, g_2)}
& = & (n^{\rho(g_1)} \eta_1 ,v^{g_1})^{(\eta_2, g_2)} \\
& = & (n^{\rho(g_1) \rho(g_2)} \eta_1^{\rho(g_2)} \eta_2 , v^{g_1 g_2}) \\
& = & (n,v)^{(\eta_1^{\rho(g_2)} \eta_2, g_1 g_2)}
\end{eqnarray*}
We must check that it preserves the adjacency structure of $\widetilde{\Gamma}$.
Observe that, by the compatibility of $\rho$ and $\phi$,
if $n_2 = n_1 \phi(v_1,v_2)$ then:
\begin{eqnarray*}
{n_2}^{\rho(g)} \phi({v_1}^g, {v_2}^g)
& = & {n_2}^{\rho(g)} \phi(v_1, v_2)^{\phi(g)} \\
& = & (n_2 \phi(v_1, v_2))^{\phi(g)} \\
& = & {n_1}^{\phi(g)}
\end{eqnarray*}
for any $(\eta,g) \in N \rtimes_\rho G$.
Also, since $G$ is a group of automorphisms of $\Gamma$, we know that
$(v_1,v_2) \in A \Gamma$ if and only if $({v_1}^g,{v_2}^g) \in A \Gamma$.
Thus the action defined above does indeed preserves adjacency.

\hfill \break
For any $(n_1,v_1), (n_2,v_2) \in V \widetilde{\Gamma}$,
by the $G$-symmetry of $\Gamma$ we can find a $g \in G$
such that ${v_1}^g = {v_2}$.
Let $\eta = {n_1}^{\rho(g)^{-1}} n_2$.
We have:
\begin{eqnarray*}
(n_1,v_1)^{(\eta,g)} & = & ({n_1}^{\rho(g)} \eta,{v_1}^g) \\
                    & = & ( {n_1}^{\rho(g)}, {n_1}^{\rho(g)^{-1}} n_2,v_2) \\
		    & = & (n_2,v_2).
\end{eqnarray*}
Thus the action is vertex transitive.

\hfill \break
Suppose that $(n_1,v_1)$ and $(n_2,v_2)$ are both adjacent to $(n,g)$,
so $n_1 = n \phi(v,v_1)$ and $n_2 = n \phi(v,v_2)$.
By the $G$-symmetry of $\Gamma$ we can find a $g \in G_v$ such that
$v_1^g = v_2$.
Let $\eta = (n^{-1})^{\rho(g)}n$. We have:
\begin{eqnarray*}
(n,v)^{(\eta,g)} & = & (n^{\rho(g)} \eta, v^g) \\
                 & = & (n^{\rho(g)} (n^{-1})^{\rho(g)}n, v) \\
		 & = & (n,v)
\end{eqnarray*}
and:
\begin{eqnarray*}
(n_1,v_1)^{(\eta,g)} & = & ({n_1}^{\rho(g)} \eta, {v_1}^g) \\
                 & = & ({n_1}^{\rho(g)} (n^{-1})^{\rho(g)}n, v_2) \\
                  & = & ({(n_1 n^{-1})}^{\rho(g)} n, v_2) \\
                   & = & (\phi(v,v_1)^{\rho(g)} n, v_2) \\
                    & = & (\phi(v,v_2) n, v_2) \\
                     & = & (n_2 n^{-1} n, v_2) \\
                      & = & (n_2, v_2) \\
\end{eqnarray*}
Thus the action is locally transitive. The result follows.	      
\end{proof}							
\end{proposition}
												 In fact, a stronger result than this holds.
												 If $\Gamma$ is $(G,s)$-arc transitive
												 then $\widetilde{\Gamma}$ is $(\widetilde{G},s)$-arc transitive,
												 where $\widetilde{G} = N \rtimes_\rho G$.
												 The proof is by induction and may be found in \cite{B1974}.
												 This construction was originally used by Conway
												 to produce an infinitely family of $5$-arc transitive graphs.

												 \hfill \break
												 For each $v \in V \Gamma$, let $B(v) = \{ (n,v) : n \in N \}$.
												 For any $(\eta, g) \in \widetilde{G}$ we have:
												 \begin{eqnarray*}
												 B(v)^{(\eta,g)} & = & \{ (n,v)^{(\eta,g)} : n \in N \} \\
												                 & = & \{ (\eta n^{\rho(g)},v^g) : n \in N \} \\
														                 & = & \{ (n, v^g) : n \in N \} \\
																                 & = & B(v^g).
																		 \end{eqnarray*}
																		 Thus the partition:
																		 $\cB = \{ B(v) : v \in V \Gamma \}$
																		 is $\widetilde{G}$-invariant.

\begin{proposition}
The quotient of the Biggs cover $\widetilde{\Gamma}$
with respect to the partition $\cB$
is isomorphic to the original graph $\Gamma$
\begin{proof}
The map $\eta : V \Gamma \to V \widetilde{\Gamma}_\cB$ given by
$v \mapsto B(v)$ establishes a bijection between the vertices of $\Gamma$
and the vertices of $\widetilde{\Gamma}_\cB$.
If $(u,v)$ is an arc of $\Gamma$, then $((1,u),(\phi((u,v)),v))$
is an arc of $\widetilde{\Gamma}$, and so $(B(u), B(v))$
is an arc of $\widetilde{\Gamma}_\cB$.
If $(u,v)$ is not an arc of $\Gamma$, then for all $n_1, n_2 \in N$
$((n_1, u),(n_2, v))$ is not an arc of $\widetilde{\Gamma}$
and so $(B(u), B(v))$ is not an arc of $\widetilde{\Gamma}_\cB$.
Thus $\eta$ establishes an isomorphism between
$\Gamma$ and $\widetilde{\Gamma}_\cB$.

\end{proof}
\end{proposition}
Let $(B,C)$ be any arc of $\widetilde{\Gamma}$.
By the above $B = B(u)$ and $C = B(v)$ for some $(u,v) \in A \Gamma$.
Its not too hard to see that each $(n,u) \in B$
has a unique neighbour in $C$, namely $(n \phi(u,v),v)$.
Thus the induced bipartite graph $\Gamma[B,C]$ is a matching.
It follows immediately that the valency of $\widetilde{\Gamma}$
is the same as the valency of $\Gamma$.

\section{Group Theoretic analysis of covering graph construction}

The Bigg's covering graph construction encompasses
all pairs of graphs $\Gamma(G, H, a)$ and $\Gamma(\widetilde{G}, H, \widetilde{a})$
where $\widetilde{G}$ is a semidirect product of $N$ by $G$ for some $N$.
The local permutation groups of $\Gamma$ and $\widetilde{\Gamma}$
are the same, and $N$ acts regularly on the fibers.

%In the more general case where $\widetilde{G}$ is an extention of $N$ by $G$
%but not a split extension $\Gamma(\widetilde{G}, H, \widetilde{a})$
%should be a cover of $\Gamma(G, H, a)$ in the sense that the induced
%bipartite graph is a matching, but Bigg's covering graph construction
%cannot be used to construct  $\Gamma(\widetilde{G}, H, \widetilde{a})$
%from $\Gamma(G, H, a)$.

%It should also possible to have a graph $\widetilde{\Gamma}$
%covering $\Gamma$ such that the local actions are the same degree
%but not permutation equivalent.

\section{Subgraph Extension}

I have been thinking about the idea of a ``combinatorial''
description of the extension of a graph in terms of its quotient.
That is, a description of the adjacency structure of $\widetilde{\Gamma}$
directly in terms of the adjaceny structure of $\Gamma$.
This section is a sketch of an idea I had, that perhaps
there is a connection between the subgraph structures
of different $G$-symmetric graphs for the same $G$.

\begin{definition} [Subgraph]
Suppose that $\Gamma = (V,A)$ is a graph.
If $W$ is a subset of $V$
and $B$ is a subset of $(W \times W) \cap A$
then $\Upsilon = (W,B)$ is a {\it subgraph} of $\Gamma$.
We write $\Upsilon < \Gamma$ and allow for the possibility
that $\Upsilon$ is a directed graph.
\end{definition}

\begin{definition} [Stabilizer of a Subgraph]
Suppose that $(G, \Gamma)$ is a symmetric graph and
$\Upsilon$ is a directed subgraph. The {\it stabilizer}
of $\Upsilon$ in $(G, \Gamma)$ is $\Aut(\Upsilon) \cap G$.
It is denoted by $G_{\Upsilon}$.
\end{definition}

\begin{definition} [Subgraph Graph]
Suppose that $(G, \Gamma)$ is a symmetric graph,
$\Upsilon$ is a directed subgraph of $\Gamma$ 
and $a$ is an involution of $\Gamma$ which fixes an arc.
The subgraph graph  $\sub(\Gamma, \Upsilon, a)$
of $\Gamma$ with respect to $\Upsilon$ and $a$
is the graph with vertices  $V = \{ \Upsilon^g : g \in G \}$
and arcs $A = \{ (\Upsilon \cup \Upsilon^a)^g : g \in G \}$.
\end{definition}

\begin{proposition}
Subgraph graph's are symmetric, with point stabilizer $G_{\Upsilon}$.
\begin{proof}
Obvious.
\end{proof}
\end{proposition}

\begin{example} [cube from tetrahedron]
Let $\Gamma$ be the tetrahedron with $S4$ acting on it.
Let $\Upsilon$ be any directed three-cycle.
Let $a$ be any involution fixing an arc of $\Upsilon$,
then the subgraph graph $\sub(\Gamma,\Upsilon,a)$
is isomorphic to the cube.
\end{example}

\begin{definition} [Symmetric subgraph]
Suppose that $(G,\Gamma)$ is a symmetric graph.
The digraph $\Upsilon < \Gamma$ is a {\it symmetric subgraph}
of $\Gamma$ if $(G_\Upsilon, \Upsilon)$ is a symmetric graph.
\end{definition}

\begin{question}
If $(G, \Gamma)$ is a symmetric graph, for which subgroups $K$ of $G$
do their exist a subgraphs of $\Upsilon$ of $\Gamma$ such that
$G_\Upsilon = K$ ?
For which subgroups $K$ of $G$
do their exist a {\it symmetric subgraph} $\Upsilon$ of $\Gamma$ such that
$G_\Upsilon = K$ ?
\end{question}

\begin{question}
Which extensions of a given $G$-symmetric graph $\Gamma$
are isomorphic to subgraph graphs of $\Gamma$ ?
Which are isomorphic to subgraph graphs of $\Gamma$
with respect to some symmetric subgraph (digraph) of $\Gamma$ ?
\end{question}